\documentclass[12pt,a4paper]{amsart}%
\usepackage{imakeidx}
\usepackage{amsmath}
\usepackage{fullpage}
\usepackage{comment}
\usepackage{MyCommA}
\usepackage[final]{pdfpages}
\usepackage{tikz-cd}

\makeatletter
\newcommand{\Rami@sectionsymbol}[1]{%
  \crefformat{#1}{##2\S##1##3}%
  \Crefformat{#1}{##2\S##1##3}%
  \crefrangeformat{#1}{\S\S##3##1##4--##5##2##6}%
  \Crefrangeformat{#1}{\S\S##3##1##4--##5##2##6}%
  \crefmultiformat{#1}{\S\S##2##1##3}{ and~##2##1##3}{, ##2##1##3}{ and~##2##1##3}%
  \Crefmultiformat{#1}{\S\S##2##1##3}{ and~##2##1##3}{, ##2##1##3}{ and~##2##1##3}%
  \crefrangemultiformat{#1}{\S\S##3##1##4--##5##2##6}{ and~##3##1##4--##5##2##6}{, ##3##1##4--##5##2##6}{ and~##3##1##4--##5##2##6}%
  \Crefrangemultiformat{#1}{\S\S##3##1##4--##5##2##6}{ and~##3##1##4--##5##2##6}{, ##3##1##4--##5##2##6}{ and~##3##1##4--##5##2##6}%
}
\Rami@sectionsymbol{section}
\Rami@sectionsymbol{subsection}
\Rami@sectionsymbol{subsubsection}
\makeatother

\newcommand{\Rami}[1]{{{#1}}}
\newcommand{\RamiA}[1]{{{#1}}}
\newcommand{\claude}[1]{{{#1}}}
\newcommand{\claudeA}[1]{{{#1}}}
\newcommand{\RamiB}[1]{{{#1}}}
\newcommand{\RamiC}[1]{{{#1}}}
\newcommand{\RamiD}[1]{{{#1}}}
\newcommand{\RamiE}[1]{{{#1}}}
\newcommand{\RamiF}[1]{{{#1}}}
\newcommand{\EitanD}[1]{{{#1}}}
\newcommand{\EitanC}[1]{{{#1}}}

\newcommand{\DimaA}[1]{{{#1}}}
\newcommand{\DimaB}[1]{{{#1}}}

\newcommand{\NextVer}[1]{}

\newcommand{\noleft}{\left.\kern-\nulldelimiterspace}

\begin{document}
	
	\author[Aizenbud]{Avraham Aizenbud}
	\address{Avraham Aizenbud,
		Faculty of Mathematical Sciences,
		Weizmann Institute of Science,
		76100
		Rehovot, Israel}
	\email{aizenr@gmail.com}
	\urladdr{https://www.wisdom.weizmann.ac.il/~aizenr/}
	
	\author[Gourevitch]{Dmitry Gourevitch}
	\address{Dmitry Gourevitch,
		Faculty of Mathematical Sciences,
		Weizmann Institute of Science,
		76100
		Rehovot, Israel}
	\email{dimagur@weizmann.ac.il}
	\urladdr{https://www.wisdom.weizmann.ac.il/~dimagur/}

\author[Kazhdan]{David Kazhdan}
	\address{David Kazhdan,
    Einstein Institute of Mathematics, Edmond J. Safra Campus, Givaat Ram The
Hebrew University of Jerusalem, Jerusalem, 91904, Israel}
	\email{david.kazhdan@mail.huji.ac.il}
	\urladdr{https://math.huji.ac.il/~kazhdan/}
	
\author[Sayag]{Eitan Sayag}
	\address{Eitan Sayag,
    Department of Mathematics, Ben Gurion University of the Negev, P.O.B. 653,
Be’er Sheva 84105, ISRAEL}
	\email{eitan.sayag@gmail.com}
	\urladdr{www.math.bgu.ac.il/~sayage}

	\date{\today}
	
	\keywords{Reductive $p$-adic group, tempered representation, character, matrix coefficient, effective bounds, Kazhdan -- Yom-Din formula}
	\subjclass{22E50, 20G25, 22E35}

	%
	%
	%
	%
	%
	%
	%
	%


    

\title{Effective and explicit version of the Kazhdan -- Yom-Din formula for characters of tempered representations}
\maketitle
\begin{abstract} 
We \RamiD{provide} an \EitanC{effective} 
\RamiD{expression} for the character of a tempered representation in terms of its matrix coefficient. This formula is an explicit version of 
a formula proven in \cite{KYtempered} which 
generalizes Harish-Chandra's formula for characters of discrete series representations, see \cite[Theorem 9]{HC_VD}.
\end{abstract}

\tableofcontents 
\section{Introduction}
We \RamiD{provide an effective expression} for the character of a tempered representation in terms of its matrix coefficient.  This formula is an explicit version of 
a formula proven in \cite{KYtempered} which 
generalizes Harish-Chandra's formula for character\DimaB{s} of discrete series representations, see \cite[Theorem 9]{HC_VD}.

Throughout the paper $\mdef{\bfG}$ is a connected reductive group defined over a local \RamiD{non-Archimedean} field $\mdef{F}$. Let $\mdef{G}:=\bfG(F)$. \EitanC{Let $\cZ(G)$ be the center of $G$ and $G^{ad}:=G/\cZ(G)$ be the adjoint group.}

\subsection{Kazhdan 
-- Yom-Din formula for characters}
Let us briefly recall a formula for characters of tempered representations from \cite{KYtempered}.
Let $\pi$ be an irreducible (complex, \RamiE{smooth}) representation of $G$.
Let $G^{ad}_n\subset G^{ad}$ be a ``ball'' of logarithmic radius $n$ inside the adjoint group $G^{ad}$, in an  appropriate sense. See \Cref{not:settings}\eqref{it:Gn} below.

\begin{theorem}[{\cite[Proposition~1.10]{KYtempered}}, see \Cref{sec:KY} \EitanC{for} more details]\label{thm:KY}
Assume that  $\pi$ is tempered. Then there exist $d\in\Z_{\geq 0}$ and $e\in \R_{>0}$ such that for any 
matrix coefficient  $m$ of $\pi$, 
we have 

$$m(1)\chi_\pi=e\lim_{n\to\infty}\frac{1}{n^{d}}\int_{G^{ad}_n}(Ad(g)m)\,dg
$$
in the space of generalized functions,
where 
$dg$ is an appropriate Haar measure on the adjoint group of $G$.
\end{theorem}
\subsection{An explicit version of Kazhdan -- Yom-Din formula}
For $m\in C^\infty(G)$ let $$\cA_n(m):= \int_{G^{ad}_n}(Ad(g)m)\,dg$$
We prove the following explicit version of \Cref{thm:KY}:
\begin{theorem*}[Informal version of \Cref{thm:recur}]
Let $\pi$ be an irreducible (complex) representation of $G$.
There exist coefficients $\alpha_1,\dots,\alpha_r\EitanC{\in \C}$ \EitanD{such that}
for any 
    \begin{itemize}
        \item matrix coefficient $m$ of $\pi$
        \item  smooth  \DimaA{compactly supported} measure $f\in \cH(G)$, 
    \end{itemize}
    the function 
    $\phi:\Z_{>0}\to \C$ defined by 
    $$\phi(k)=\sum_{i=1}^r  \alpha_i \langle \cA_{k+i}(m),f \rangle$$
    is constant whenever $k$ is large enough.

Moreover,
if $\pi$ is tempered we can choose $\alpha_i$ such that for large enough $k$ we will have 
$$\phi(k)=m(1)\langle \chi_\pi,f\rangle. $$ 
\end{theorem*}

We prove an effective version of this result, i.e. we show that the $k$ where $\phi(k)$ becomes constant depends linearly on the level of $m$ and $f$. In order to formulate \RamiD{this version} properly we will need to introduce some notation.

\subsection{Notation}\label{ssec:intro.not}
For a group  (or an algebraic group) $H$ denote by  $\mdef{\cZ(H)}$ its center.    For  a smooth (complex) representation $\pi$ of $G$ and $v\in \pi$,  $w\in \tilde\pi$, we denote by $\mdef{m_{v,w}}:G\to \C$ the \EitanC{corresponding} matrix coefficient.

In order to formulate our results we have to choose some additional data related to $\bfG$.\footnote{\RamiA{Though the dependence of our results on these choices is rather weak, the choices cannot be made completely arbitrarily.}} 
In \Cref{ssec:not} we specify how we make these choices. For now we will just 
\RamiF{list them.}
We choose:
\begin{itemize}
\item  A maximal compact subgroup $G_0<G$.
\item  A basis of the topology of $G$ at $1$ consisting of open compact subgroups $G_{-n}<G_{-(n-1)}<G_0$ for $n\in \Z_{> 0}$.
    \item Open compact subsets $G_n^{ad}\subset G^{ad}_{n+1}\subset\mdef{G^{ad}}:=G/\cZ(G)$ \EitanC{that} cover $G^{ad}$.
\end{itemize}
We also \EitanC{introduce the following notation}:
\begin{itemize}
    \item 
    For an open  compact subgroup $K<G$ denote by 
$\mdef{\cH_{K}(G)}$ the \tdef{Hecke algebra} consisting of compactly supported measures on $G$ that are bi-$K$-invariant.
    \item Let $\mdef{\mu_{G^{ad}}}$ be the Haar measure on $G^{ad}$ normalized on  the image of $G_0$ in $G^{ad}$.
    \item For $f\in C^{\infty}(G)$  and an integer $n\in \Z_{>0}$ denote by $\mdef{\cA_n(f)}\in C^{\infty}(G)$ the function given by $$(\cA_n(f))(x):=\int_{g\in G^{ad}_n} f(gxg^{-1}) \mu_{G^{ad}}$$
\end{itemize}

\subsection{Main results}
\begin{introtheorem}[\Cref{ssec:KY}]\label{thm:recur}
Let $\pi$ be an irreducible  representation. 
Then there exist\footnote{\RamiA{In fact, the integer $r$ can be bounded in a way that depends only on the Weyl group of $G$  (and not on the representation $\pi$), cf. \cite[Theorem C]{AGKS_eff_rep}.}} coefficients $\alpha_1,\dots,\alpha_r\EitanC{\in \C}$ \RamiC{(not all $0$)}  and a constant $\RamiD{c}\in \Z_{>0}$ s.t. 
for any 
    \begin{itemize}
        \item $n\in \Z_{>0}$
    \item $v \in \pi^{G_{-n}}$,  $w\in \tilde\pi^{G_{-n}}$
        \item $f\in \cH_{G_{-n}}(G)$,
    \end{itemize}
    the function 
    $\phi:\Z_{>0}\to \C$ defined by 
    $$\phi(k)=\sum_{i=1}^r  \alpha_i \langle \cA_{k+i}(m_{v,w}),f \rangle$$
    is constant whenever $k> \RamiD{c}
    \cdot n$.
Moreover,
if $\pi$ is tempered then we can choose $r$ and the \RamiF{$\{\alpha_i\}_{i=1}^r$} such that for \EitanD{$k>c
    \cdot n$}  we have: 
$$\phi(k)=\langle \chi_\pi,f\rangle \langle v,w \rangle.$$
\end{introtheorem}
\subsection{Motivation}
\RamiC{This paper is part of 
our project of extending Harish-Chandra's integrability theorem from \EitanD{the case of reductive groups over non-archimedean fields of characteristic zero} to reductive groups over local fields of positive characteristic (\cite{AGKS2,AAGKS,AGKS_jet,AGKS_orb,AGKS_eff}).}

Specifically, we plan to use  \Cref{thm:recur} to generalize the results of \cite{AGKS2,AAGKS} to arbitrary irreducible representations.

\RamiA{In \cite{HC55,HC_VD} Harish-Chandra gives a formula for characters of discrete series  representations in terms of matrix coefficients. Moreover, he indicated an effective version for cuspidal representations \RamiF{of $\GL_n(F)$}. We used the effectivity of this formula 
in \cite{AGKS2}, in order to prove local integrability of characters of cuspidal representations}. \RamiA{The fact that Harish-Chandra's formula for the characters of discrete series representations is not effective} \EitanC{makes the arguments of \cite{AGKS2} inapplicable for discrete series representations}.

We plan to use \Cref{thm:recur} in order to generalize the results of 
\cite{AGKS2,AAGKS} \RamiA{concerning local integrability of characters from cuspidal representations} to tempered representations, see \cite{AGKS_int}.  This should imply \RamiF{the generalization of} the results {of 
\cite{AGKS2,AAGKS}} \RamiF{to} arbitrary \EitanD{irreducible} representations. 

\subsection{Related results}
\RamiA{In the literature, there are several formulas for the character of an irreducible representation in terms of {its} matrix coefficients, in various generalities. To the best of our knowledge, the most  general one is given 
in \cite{KYtempered}, and it is valid for characters of tempered representations.}

One can organize the formulas for the character \EitanC{of an irreducible representation} in terms of \EitanC{its} matrix coefficients into two types:
\begin{enumerate}
    \item formulas for
the character as a distribution.
\item
formulas for the value of the character at nice points.
\end{enumerate}
This paper focuses on the first type.
\EitanC{In \cite[Theorem 9]{HC_VD}, Harish-Chandra provided a formula of \EitanC{the first} type for cuspidal representations.}  In characteristic \EitanD{zero} one can use this formula (together with an additional result of Harish-Chandra concerning the stabilization of the averaging of matrix coefficients, see \cite[page 101]{HC_VD}) to obtain a formula of the second type. Note that such a formula will be more effective than the one in \cite[Theorem 9]{HC_VD}. This approach was extended to positive characteristic in \cite{AGKS_orb}. 

Another formula for characters of the second type in the generality of 
 elliptic tempered representations was developed in \cite{Art76, Art87, Art93} by a different approach. However, in order to derive information from this formula for the character as a distribution, one needs the local integrability of characters, and hence this approach is not valid (yet) in positive characteristic and cannot be used for proving local integrability of characters.

Another related work \EitanC{is} \cite{BezKaz} where the case of $G=\SL_2$ is studied.

As mentioned above \cite{KYtempered} proves a version of Harish-Chandra's formula for tempered representations, and our proof of \Cref{thm:recur} uses that formula.

In \cite{AGKS_eff_rep} we provided \EitanD{the} infrastructure for proving effective results of the type of \Cref{thm:recur}. We use this infrastructure here.
 
\subsection{Idea of the proof of \Cref{thm:recur}}\label{sub:idea}
Following 
\cite{KYtempered, JL, HC_VD} 
we use the Cartan decomposition in order to express $\langle \cA_{n}(m_{v,w}),f \rangle$ as a sum of integrals over $G_0$-double cosets. 

\EitanC{These double cosets are indexed by the integral part $\Lambda^+$ of the Weyl chamber. 
\RamiA{For simplicity, assume that $\Lambda^+$ is the collection of the dominant co-characters of a maximal split torus in $G$. \RamiF{So} we have a map $\lambda \mapsto \varpi^\lambda$ from $\Lambda^+$ to $G$.}

For each $\lambda\in \Lambda^+$ denote the corresponding integral by $\cS^1(\lambda)$.} 
\RamiA{Namely:
$$\cS^1(\lambda):=\int_{G_0 \varpi^\lambda G_0} \langle Ad(h)(m_{v,w}),f\rangle \,\mu_{G^{ad}}.$$
}

\EitanC{Along each ray in $\Lambda^+$,} one can present the function $\cS^1$ as a combination of restrictions of matrix coefficients to $\Lambda^+$ \EitanC{(such an equality is valid far enough in the ray).} 
\RamiA{More explicitly,
we introduce a function $\cS:\Lambda^+\to \C$ given by $$\cS(\lambda)=\langle Ad(\varpi^\lambda)(m_{v,w}),f\rangle,$$
and, far enough along each ray, we  express $\cS^1$ as a combination of functions like $\cS$ \EitanD{(with different parameters $v,w,f$)}.
}

\EitanC{
\RamiA{The function $\cS$} satisfies an effective recursion relation along each of the the directions in $\Lambda^+$. This \RamiA{follows from} \cite[Theorem C]{AGKS_eff_rep}, see \Cref{cor:S.adm} for more details.}

\EitanC{We deduce the first part of \Cref{thm:recur} (before the ``\EitanD{M}oreover'') in two steps:
\begin{itemize}
    \item Prove that $\cS^1$ satisfies an effective recursion relation along each of the directions in $\Lambda^+$.
    \item Deduce that the sequence of partial sums of $\cS^1$ over ``balls" in $\Lambda^+$ satisfies an effective recursion relation.
\end{itemize}
To prove these recursions, we introduce and study the notion of an $(\cR,n)$-\EitanD{\it admissible} function on an integral cone $L$. This notion informally means that the function satisfies an effective recursion relation along each of the directions in the cone.
\RamiA{The parameters $\cR,n$ provide a quantitative  description of the recursions,} see  \Cref{def:adm}.

The key tool to establish the above steps is the
following  
\EitanD{result on \it{admissible} functions} : 
}
\EitanC{
\begin{introtheorem}[\Cref{cor:cone.prod}, \Cref{cor:acc_pushforward}, \Cref{lem:fin.index}]\label{thm:cones}
    Let $L\subset \Z^N$ be a (strict) integral cone. 
     Then:
     

    \begin{enumerate}
        \item         
     Suppose $f,g\in \C^L$ are two  admissible functions. Then $fg$ is admissible.
        \item          
        Suppose $f\in \C^L$ is an admissible  function. Let $e:L\to \Z_{\geq 0}$ be a monoid homomorphism with trivial kernel. Then the function $$k\mapsto \sum_{e(x)\leq k} f(x)$$ is admissible.
 \item Let   $f\in \C^L$ be a function. Let $L'\subset L$ be a subcone of  a finite index. 
 \DimaB{If} $f$ is admissible on every coset of $L'$, then it is admissible on $L$.
    \end{enumerate}
\end{introtheorem}
\RamiF{T}he rigorous formulation of this theorem includes control on the  parameters \EitanD{$(\cR,n)$ of the admissibility of the various functions involved.}
}

The ``\EitanD{M}oreover'' \EitanC{part of \Cref{thm:recur}} follows from the main part using \Cref{thm:KY}  and the elementary fact that a limit of a sequence that satisfies a recursion relation can be expressed in an effective way as a \RamiE{fixed} finite combination of elements of the sequence (see \Cref{lem:rec}).

\subsection{Structure of the paper}
In \Cref{sec:prelim} we fix our conventions, we make the choices of the data listed in
\Cref{ssec:intro.not}, we introduce the notation attached to them, and we record the
basic facts on these data that are used throughout the paper. In \Cref{ssec:eff.rep} we
recall the effective recursion satisfied by the matrix coefficients of a
representation of finite length along a dominant cocharacter (see \Cref{thm:mcexp}). \EitanD{This is a result we proved in
\cite{AGKS_eff_rep}.}

In \Cref{sec:cones} we prove \Cref{thm:cones}. 

In \Cref{sec:adm.mat} we prove 
\Cref{cor:S.adm} that states that 
\RamiA{the function $\cS$ introduced in \Cref{sub:idea}} 
is $(\cR,n)$-admissible for some fixed $\cR$ (depending on $\pi$) and for $n$ that depends linearly on the \EitanC{level}\RamiF{s} of \RamiA{$v,w$} and $f$. This is based on the results of \cite{AGKS_eff_rep}.

In \Cref{ssec:KY} we prove \Cref{thm:recur}. 
\EitanC{We analyze the function $\cS^1$ introduced in \Cref{sub:idea}. In order to present it as a combination of \RamiA{some variations of the} function $\cS$, we need to analyze the following:
\begin{itemize}
    \item The volume of a $G_{-n}$-double coset inside a given  $G_0$-double coset.
    \item \RamiD{The decomposition of a given $G_0$-double coset into} $G_{-n}$-double cosets.
\end{itemize}

 In  \Cref{ssec:vol} we show  that the above volume is an admissible function of the $G_0$-double coset.
  In \Cref{ssec:stab} we show that the stabilizer of a $G_{-n}$-double coset \EitanD{with respect to} the action of $G_{0}\times G_{0}$
stabilizes along a ray in the Weyl chamber  after a number of
steps proportional to $n$. This provides \RamiD{a certain stabilization property} of the collection of $G_{-n}$-double cosets inside a given $G_0$-double coset.
}
In
\Cref{ssec:avg} we combine this with \Cref{thm:mcexp} \EitanD{(\it{matrix coefficients are admissible functions})} and with \Cref{cor:cone.prod} \EitanD{(\it{product of admissible functions is admissible})} and
prove that $\cS^1$ is an admissible function when restricted to cosets of split elements (\Cref{prop:adm.int}). In \S \ref{ssec:avg.G} we deduce from this that $\cS^1$ is admissible.
In \Cref{ssec:lim} we record the
elementary fact that a limit of a sequence satisfying a recursion relation is an effective
finite combination of the elements of the sequence (\Cref{lem:rec}). In
\Cref{ssec:pf.recur} we deduce \Cref{thm:recur}.

In \Cref{sec:KY} we \EitanD{provide the necessary \claude{adaptations} to deduce \Cref{thm:KY} from the proof of \cite[Proposition~1.10]{KYtempered}. }
\subsection{Acknowledgments}
During the preparation of this paper, A.A., D.G. and E.S. were partially supported by the ISF grant no. 1781/23. 
D.K. was partially supported by an ERC grant 101142781.

\section{Preliminaries}\label{sec:prelim}
\subsection{Conventions}
Here we fix general notation and conventions that we will use throughout the paper.

\begin{enumerate}
    \item By a \tdef{variety} we mean a reduced scheme of finite type over $F$. Similarly for an algebraic group.

    \item We will usually denote algebraic varieties by bold face letters (such as $\bfX$) and the spaces of their $F$-points by the corresponding usual face letters (such as $X:=\bfX(F)$). We use the same conventions when we want to interpret vector spaces as algebraic varieties.

    \item We will use the same letter to denote a morphism between algebraic varieties and the corresponding map between the sets of their $F$-points.
\item By a representation of an $l$-group we will mean a complex-valued smooth representation.
    \item {We will use the symbol \tdef{$<$} to denote the (not necessarily proper) containment relation for groups}.
\item \Rami{For a group $\Gamma$ and subsets $X_1,\dots,X_N\subset \Gamma$ we denote by $\langle X_1,\dots,X_N \rangle$ the subgroup of $\Gamma$ generated by \EitanD{$X_1 \cup \dots \cup X_N$.} We will use a similar notation for algebraic groups and other structures.}
\item \Rami{For a group $\Gamma$ acting on a space $X$ and $x\in X$ we denote by $\Gamma_x$  the stabilizer of $x$ \EitanD{in $\Gamma$}. We will also use a similar notation for algebraic groups.}
\item We will use the words parabolic/Levi to describe subgroups both of algebraic reductive groups and of the groups of their points.

\item 

Several statements in this paper will concern the existence of objects that satisfy some conditions. Sometimes we will want to fix such an object. In these cases 
 we assign a distinct name to the corresponding object inside the formulation of the statement. It is implied that after each such statement we fix such an object and we will refer to it later by this name. Nothing significant will depend on the choices of these objects.

\end{enumerate}

\subsection{Notation and choices of the data}\label{ssec:not}
Let us fix the choices described in \Cref{ssec:intro.not} and introduce some notation related to them.
\begin{notation}\label{not:settings}$ $
    \begin{enumerate}
         \item Let $\mdef{\ell}$ be the cardinality of the residue field of $F$.
        \item Fix a uniformizer $\mdef{\varpi}\in F $.
        \item Fix a  minimal parabolic $\mdef{\bfB}\subset \bfG$.
        \item Fix a  Levi subgroup $\mdef{\bfT}<\bfB$.
        \item Fix a  maximal split torus $\mdef{\bfA_\bfT}<\bfT$.       
        \item Denote by $\mdef{\cB\cT(G)}$ the (extended) Bruhat-Tits building of $G$ \RamiF{and by $\mdef{\cB\cT_{red}(G)}$ the reduced Bruhat-Tits building of $G$}.
        \item\label{it:vertex} Fix a special vertex $\mdef{\fx}\in  \cB\cT(G)$ that lies in the apartment that corresponds to $\bfA_{\bfT}$. Such a vertex exists by \cite[Proposition 1.3.43]{KP}.
        \item We set $\mdef{G_0}:=G_\fx.$ 
        \item \RamiA{For $n\in \Z_{>0}$ w}e set $\mdef{G_{-n}}:=\mdef{G_{\fx,n}}$ to be the \tdef{Moy-Prasad group}. (See \cite{MP96} or \cite[Definition 13.2.1]{KP}.)        
        \item\label{it:0H} For a connected reductive group $\bfH$  with $H=\bfH(F)$ denote by $\mdef[^0G,\, ^0H,\, ^0M,\, ^0T,\, ^0A]{^0H}$ the 
        intersection of the kernels of all unramified characters of $H$. See e.g. \cite[V.2.3]{Ren}. Every compact subgroup of $H$ is contained in $^0H$ (\cite[V.2.3]{Ren}). In particular, if $H$ has no proper parabolics, so that $^0H$ is compact, then $^0H$ is the unique maximal compact subgroup of $H$.        
        \item For a connected reductive group  $\bfH$ denote $\mdef{\Lambda(\bfH)}:=\mdef{\Lambda(H)}:=H/^0H$ where $H=\bfH(F)$. For non-connected $\bfH$ we denote $\Lambda(\bfH)=\Lambda(\bfH^\circ)$, where $\bfH^\circ$ is the identity component of $\bfH$.        
                \item\label{it:Phi} Denote by $\mdef{\Phi}$ the set of non-divisible roots of the root system defined by $(\bfG,\bfA_\bfT)$, i.e. those roots $\beta$ of $\bfA_\bfT$ in $\Lie(\bfG)$ for which $\beta/2$ is not a root. It is itself a root system. We consider its elements as $\Q$-valued functionals on $\Lambda(T)$. See \cite[V.3.7, V.3.17, V.2.6]{Ren}.         
        \item A \tdef{standard parabolic subgroup} of $\bfG$ is a parabolic that contains $\bfB$. 
        A  \tdef{standard Levi subgroup} of $\bfG$ is a Levi subgroup  of a standard parabolic subgroup of $\bfG$ that contains $\bfT$. We apply analogous terminology for parabolic and Levi subgroups of $G$.
\item\label{it:opp} For a standard parabolic $\bfP<\bfG$ we denote by $\mdef{\bfP^-}$ the unique parabolic s.t. $\bfP\cap \bfP^{-}$ is a standard Levi. We call it the \tdef{opposite parabolic} of $\bfP$.

\item\label{it:Phi.P} For $\bfP<\bfG$ a standard parabolic or the opposite of one, with unipotent radical $\bfU_\bfP$, denote by $\mdef{\Phi_\bfP}\subset\Phi$ the set of $\beta\in\Phi$ occurring in $\mdef{\bfU_\bfP}$, i.e. those whose weight space for $\bfA_\bfT$ lies in $\Lie(\bfU_\bfP)$. Then $\Phi=\Phi_\bfB\sqcup\Phi_{\bfB^-}$ and $\Phi_{\bfP^-}=-\Phi_\bfP$. If $P=\bfP(F)$ we use the same notation with $\bfP$ replaced by $P$ to denote the same set.
        \item\label{it:Delta} Denote by $\mdef{\Delta}\subset\Phi_\bfB$  the collection of \tdef{simple roots} of the root system defined by $(\bfB,\bfA_{\bfT})$. See \cite[V.3.10]{Ren}.
        \item\label{it:Weyl} Let $\mdef{\cN}$ be the normalizer of $T$ in $G$ and let $\mdef{W}:=\cN/T$ be the \tdef{Weyl group} of the root system $\Phi$ (see \cite[\S V.3]{Ren}). We let $W$ act on
        $\Lambda(T)$, 
        by $\RamiE{\mdef{w_\nu}([t]):=[\nu t\nu^{-1}]}$, where $w_\nu\in W$ denotes the class of $\nu\in\cN$.
        \item\label{it:M+} Let $$\mdef{T^+}=\{t\in T|\forall \beta\in \Phi_\bfB \text{ we have } \beta(t)\geq 0  \},$$ see \cite[\S V.3.19]{Ren}.
        \item\label{it:Lam+} Denote $\mdef{\Lambda^+(T)}:= T^+/^0T$.
        \item\label{it:Lam+A} Pulling back to $A_{\bfT}$ we get $\mdef{\Lambda^+(A_{\bfT})}\Rami{:=(A_{\bfT}\cap T^+)/{^0A}_{\bfT}}$.
                \item Denote by $\mdef{\cE}:=\mdef{\cE_{A_{\bfT}}}:\Lambda({A_{\bfT}})\to {A_{\bfT}}$ the natural map given by the choice of $\varpi$. We also denote $\mdef{\varpi^\lambda}:=\cE_{A_{\bfT}}(\lambda)$
        
\item\label{it:N.Theta}        Choose $\mdef{\fN}\in\Z_{>0}$ such that
$\fN\beta(\lambda)\in\Z$ for every $\lambda\in\Lambda(T)$ and every $\beta\in\Phi$, the elements
of $\Phi$ being viewed as $\Q$-valued functionals on $\Lambda(T)$ as in \eqref{it:Phi}. Such an
$\fN$ exists since $\Phi$ is finite and $\Lambda(T)$ is finitely generated.
Denote
$$
\mdef{\Theta}=\fN\sum_{\beta\in\Phi_\bfB}\beta.
$$
        \item Let $\mdef{G_0^{ad}}$ be the image of $G_0$ in $G^{ad}$. 
        \item\label{it:cartan} $T^+\to G$ induces a bijection 
        $\Lambda^+(T)\to G_0 \backslash G/ G_0 $ (\tdef{Cartan decomposition}, \cite[\S3.3.3]{Tit79}
        ).
        \item\label{it:ZG} \RamiE{Note that $\cZ(G)\subset T^+$. Denote the  image of $\cZ(G)$ in $\Lambda^+(T)$ by $\Lambda_{\cZ(G)}(T)$.}
        \item\label{it:ov1} The map in \eqref{it:cartan} induces a bijection $\Lambda^+(T)/\RamiE{\Lambda_{\cZ(G)}(T)}\to G_0^{ad} \backslash G^{ad}/ G_0^{ad}$. 
         Note also that $\Theta$ is trivial on $\claude{\Lambda_{\cZ(G)}(T)}$. So it gives a function on $\Lambda^+(T)/\claude{\Lambda_{\cZ(G)}(T)}$ and hence on $G_0^{ad} \backslash G^{ad}/ G_0^{ad}$ and on $G^{ad}$. We denote this function by $\mdef{ov_{G^{ad}}}$
        \item\label{it:ov2} Let $\mdef{\bfG^{ab}}:=\bfG/[\bfG,\bfG]$ and let  $\mdef{G^{ab}}:=\bfG^{ab}(F)$.
        \item\label{it:Gn} Let $\mdef{G^{ad}_n}:=ov_{G^{ad}}^{-1}([0,n])$.
     \end{enumerate}
\end{notation}
\subsection{Basic facts on the data}

\begin{lemma}[{\cite[Theorem 4.2]{MP96}}]\label{lem:iva.decomp}
    Let $P<G$ be a standard parabolic. Let $M<P$ be its standard Levi subgroup and $U<P$ be its unipotent radical. Let $U^-$ be the unipotent radical of $P^-$. Then
    for \DimaA{every} $n>0$ we have 
    $$G_{-n}=(U^-\cap G_{-n})(M\cap G_{-n})(U\cap G_{-n}).$$
\end{lemma}

\begin{notation}\label{not:root.sub}
$ $
\begin{enumerate}
\item\label{not:root.grp}
For every $\beta\in \Phi$, let $\mdef{\bfU_\beta}$ be the corresponding \tdef{root subgroup}, i.e. the subgroup
whose Lie algebra is the sum of the weight spaces of $\bfA_\bfT$ in $\Lie(\bfG)$ whose weights
are the positive integer multiples of $\beta$. See \cite[V.3.7]{Ren}.
\Rami{Let $\mdef{U_\beta}:={\bfU_\beta(F)}$.}
\item\label{not:root.filtration}
For $\beta\in\Phi$ and $r\in\R$ let
$$
\mdef{U_{\beta,r}}:=U_{\beta,\fx,r}\subset U_\beta
$$
be the Bruhat--Tits root-filtration subgroup at the fixed special vertex
$\fx$, i.e. $U_{\beta,r}=\phi_\beta^{-1}([r,\infty])$ for the valuation $\phi$ of the root
datum of $(\bfG,\bfA_\bfT)$ corresponding to $\fx$ \RamiA{(see \cite[Construction 9.6.1]{KP})}. We use the normalization of
\cite[Definition 6.1.2]{KP}. The
family $(U_{\beta,r})_{r\in\R}$ is decreasing in $r$.
\item\label{not:root.div}
\RamiD{Let $\beta\in\Phi$ be such that $2\beta$ is a root of $(\bfG,\bfA_\bfT)$, let
$\mdef{U_{2\beta}}<U_\beta$ be the root subgroup of $2\beta$, and for $r\in\R$ let
$$
\mdef{U_{2\beta,r}}:=U_{2\beta,\fx,r}=\phi_{2\beta}^{-1}([r,\infty])\subset U_{2\beta}
$$
be the corresponding root-filtration subgroup at $\fx$. 
If $2\beta$ is not a root we set $U_{2\beta}:=\{1\}$ and $U_{2\beta,r}:=\{1\}$, as in
\cite[Remark 7.3.4]{KP}. By axiom V\,4 of \cite[Definition 6.1.2]{KP} we have
$\phi_{2\beta}=2\phi_\beta$ on $U_{2\beta}$, and consequently
$$
U_{2\beta,2r}=U_{2\beta}\cap U_{\beta,r}\qquad(r\in\R).
$$}

\end{enumerate}
\end{notation}

\begin{lemma}[{\cite[Definition 6.1.2, Axiom V\,6]{KP}}]\label{lem:root.conj}
$T$ normalizes each $U_\beta$ and shifts its filtration by the value of $\beta$:
$$
tU_{\beta,r}t^{-1}=U_{\beta,\,r+\beta(t)}\qquad(t\in T,\ \beta\in\Phi,\ r\in\R),
$$
where $\beta(t)$ is the value at the class of $t$ in $\Lambda(T)$ of $\beta$ viewed as a
$\Q$-valued functional as in \Cref{not:settings}\eqref{it:Phi}.
\end{lemma}

\begin{lemma}[{\cite[Lemma A.4.\RamiA{6} \RamiD{and its proof,} \RamiF{and Lemma A.4.4}]{AGKS_eff_rep},\cite[7.3.12(3)]{KP}}]
\label{lem:U.filtration}
Let \Rami{$\bfQ<\bfG$} be a
standard parabolic or the opposite of one, let \Rami{$\bfU_\bfQ$} be its unipotent radical, and
fix an \RamiD{arbitrary} ordering of $\Phi_{\Rami{\bfQ}}$. 
\Rami{Let $Q:=\bfQ(F)$ and $U_Q:=\bfU_\bfQ(F)$.}
Then, with the notation of \Cref{not:root.sub}:
\RamiD{
\begin{enumerate}
\item\label{lem:U.filtration:1} 
{For every $n\in\Z_{\geq 0}$ the multiplication map
$$
\prod_{\beta\in\Phi_Q}U_{\beta,n}\longrightarrow U_Q
$$
is injective, and its image satisfies
$$
U_Q\cap G_{-2n}\subset\prod_{\beta\in\Phi_Q}U_{\beta,n}\subset U_Q\cap G_{-n}.
$$
}
\item\label{lem:U.filtration:3}
{Let $s\in T$ satisfy $\beta(s)\geq0$ for every $\beta\in\Phi_Q$. Then for every
$n\in\Z_{\geq0}$
$$
s\,(U_Q\cap G_{-n})\,s^{-1}\subset U_Q\cap G_{-n}.
$$}
\item\label{lem:U.filtration:prod} 
For every $n\in\Z_{\geq 0}$ the multiplication map
$$
\prod_{\beta\in\Phi_Q}(U_{\beta,n}U_{2\beta,n})\longrightarrow (U_Q \cap G_{-n})
$$
is a bijection. 
\item\label{lem:U.filtration:4} \claude{The multiplication map
$$
\prod_{\beta\in\Phi_Q}U_{\beta}\longrightarrow U_Q
$$
is injective.}
\end{enumerate}
}
\end{lemma}
\begin{lemma}\label{lem:parahoric}$ $
\begin{enumerate}
\item\label{lem:parahoric:norm} For every $n\in\Z_{\geq 0}$ the group $G_{-n}$ is normal in $G_0$.
\item\label{lem:parahoric:T} $T\cap G_0={}^0T$.
\item\label{lem:parahoric:fix} For every $t\in T$
$$
G_0\cap tG_0t^{-1}=\bigl\langle\, \cN\cap G_0\cap tG_0t^{-1},\ \ \Rami{\bigcup_{\beta\in\Phi}} U_{\beta,\max(0,\beta(t))}\, \bigr\rangle ,
$$
where $\beta(t)$ is as in \Cref{lem:root.conj}.
\end{enumerate}
\end{lemma}
\begin{proof}
\claude{Recall that $G_0=G_\fx$ is the stabilizer of $\fx$ in $G$ acting on the extended building. By \cite[Lemma 4.3.3]{KP}\footnote{\RamiF{In \cite{KP} our field $F$ is denoted by $k$ and our algebraic group $\bfG$ is denoted by $G$. So $G(k)$ in \cite{KP} means our $G$.}} it is the group $\RamiF{G}^1_\fx$ of \cite[(7.7.1)]{KP} attached to the image of $\fx$ in \claudeA{the reduced building $\cB\cT_{red}(G)$}.  \RamiF{Here} $\RamiF{G}^1={}^0G$. As $\fx$ is special, $\RamiF{G}^1_\fx=\RamiF{G}^b_\fx=\RamiF{G}^\dagger_\fx$} (\RamiA{see \cite[Lemma 4.3.3, Proposition~7.7.11]{KP}}).

\claude{\eqref{lem:parahoric:norm} For $n=0$ the assertion \RamiF{is clear}. Let $n>0$. By \cite[Proposition~13.2.5(2)]{KP} we have $gG_{\fx,n}g^{-1}=G_{g\fx,n}$ for every $g\in G$. \RamiF{Since} every $g\in G_0$ fixes $\fx$ \RamiF{we have} $gG_{-n}g^{-1}=G_{-n}$.}

\claude{\eqref{lem:parahoric:T} 
\RamiF{Follows from \Cref{not:settings}\eqref{it:cartan}}.
}

\claude{\eqref{lem:parahoric:fix} The element $t$ acts on the apartment corresponding to $\bfA_\bfT$ by a translation (\cite[Axiom~4.1.4(A\,1)]{KP}), so $\Omega:=\{\fx,t\fx\}$ is a non-empty bounded subset of it. Applying \cite[Lemma 4.3.3]{KP} to $t\fx$ gives $tG_0t^{-1}=\RamiF{G}^1_{t\fx}$, whence
$$
G_0\cap tG_0t^{-1}=\RamiF{G}^1_\fx\cap \RamiF{G}^1_{t\fx}=:\RamiF{G}^1_\Omega
$$
by \cite[Lemma 7.7.4]{KP}. By the construction \cite[(7.7.1)]{KP}, $\RamiF{G}^1_\Omega$ is the product of the 
$\RamiF{\cN}^1_\Omega\RamiF{:=\cN\cap \RamiF{G}^1_\Omega}$
with the subgroup generated by the groups $U_{\beta,\Omega,0}\RamiF{:=U_{\beta,\fx,0}\cap U_{\beta,t\fx,0}}$ \RamiF{where} $\beta$ \RamiF{runs} over the roots of $(\bfG,\bfA_\bfT)$ (\cite[Definition~7.4.1(2)]{KP}). For $\beta\in\Phi$ we have
$$
U_{\beta,\Omega,0}=U_{\beta,\fx,0}\cap U_{\beta,t\fx,0}
\overset{\text{Lem }\ref{lem:root.conj}}{=}U_{\beta,0}\cap U_{\beta,\beta(t)}
\overset{\text{Not }\ref{not:root.sub}\eqref{not:root.filtration}}{=}U_{\beta,\max(0,\beta(t))},
$$
the first equality by \cite[Definition~7.4.1(1)]{KP}.
\RamiF{Also} $U_{2\beta,\Omega,0}\subset U_{\beta,\Omega,0}$ by \Cref{not:root.sub}\eqref{not:root.div}. Finally
\RamiF{
\begin{align*}
G_0\cap tG_0t^{-1}&={G}^1_\Omega\\
&=\langle 
\cN^1_\Omega,\bigcup_{\beta\in \Phi \cup 2\Phi} U_{\beta,\Omega,0}\rangle\\
&\overset{\claude{\text{Not }\ref{not:root.sub}\eqref{not:root.div}}}{=}
\langle 
\cN \cap  {G}^1_\Omega,\bigcup_{\beta\in \Phi} U_{\beta,\Omega,0}\rangle\\
&=
\langle 
\cN \cap G_0\cap tG_0t^{-1} ,\bigcup_{\beta\in \Phi} U_{\beta,\max(0,\beta(t))}\rangle
\end{align*}
as required.
}
}
\end{proof}

\begin{lemma}\label{lem:T.norm}
\claude{For every $n\in\Z_{>0}$ the group $T\cap G_{-n}$ is normal in $T$.}
\end{lemma}
\begin{proof}
\claude{For a point $y$ of the apartment of $\cB\cT(G)$ corresponding to $\bfA_\bfT$ let $G_{y,n}$ be the Moy--Prasad group of \cite[Definition 13.2.1]{KP}, so that $G_{\fx,n}=G_{-n}$, and let $T_n<T$ be the group denoted $Z(k)_n$ in \cite[Proposition 13.2.5(3)]{KP}.}

\claude{We first claim that 
\begin{equation}\label{eq:T.Gn.Tn}
 T\cap G_{y,n}=T_n   
\end{equation}
for every such $y$. Indeed $T_n<T\cap G_{y,n}$ by \cite[Proposition 13.2.5(3)]{KP}. Conversely let $z\in T\cap G_{y,n}$. Applying \cite[Proposition 13.2.5(3)]{KP} with the system of positive roots $\Phi_\bfB$ we may write $z=u^-\tau u^+$, where $\tau\in T_n$, where $u^-$ is a product of elements of the groups $U_{\beta,y,n}$ with $\beta\in\Phi_{\bfB^-}$, and where $u^+$ is a product of elements of the groups $U_{\beta,y,n}$ with $\beta\in\Phi_\bfB$. In particular $u^-\in U^-$ and $u^+\in U^+$, the unipotent radicals of $B^-$ and of $B$. Comparing with the factorization $z=1\cdot z\cdot 1$ and using the injectivity of the multiplication map $U^-\times T\times U^+\to G$ we get $u^-=u^+=1$ and $z=\tau\in T_n$.}
\RamiE{
Now let $g\in T$. By \cite[Axiom 4.1.4(A\,1), \S4.3]{KP} the group $T$ preserves the apartment corresponding to $\bfA_\bfT$, so $g\fx$ lies in that apartment. Hence}
\RamiF{
$$g(T\cap G_{-n})g^{-1}=T\cap gG_{\fx,n}g^{-1}\overset{\text{\cite[Prop. 13.2.5(2)]{KP}}}{=}T\cap G_{g\fx,n}\overset{\eqref{eq:T.Gn.Tn}}{=}T_n\overset{\eqref{eq:T.Gn.Tn}}{=}T\cap G_{-n}.$$

}

\end{proof}

\begin{lemma}\label{lem:T.conj}
Let $X\subset T$ be a finite set. Then there is $c\in\Z_{>0}$ s.t. for any $g\in X$ and any $n\in\Z_{>0}$ we have
$$gG_{-cn}g^{-1}<G_{-n}.$$
\end{lemma}
\begin{proof}
\RamiE{WLOG \claude{assume} $\Phi\neq \emptyset$.}
By \Cref{not:settings}\eqref{it:Phi.P} we have $\Phi=\Phi_{B^-}\sqcup\Phi_{B}$. Put
$$c_\RamiD{0}:=1+\max_{g\in X}\max_{\beta\in\Phi}\lceil |\beta(g)|\rceil\in\Z_{>0},$$
and
$$c=2c_0.$$

Fix $g\in X$ and $n\in\Z_{>0}$. For every $\beta\in\Phi$ we have $c_\RamiD{0}n+\beta(g)\geq c_\RamiD{0}n-(c_\RamiD{0}-1)\geq c_\RamiD{0}n-(c_\RamiD{0}-1)n=n$, hence
\begin{equation}\label{lem:T.conj:2}
gU_{\beta,c_\RamiD{0}n}g^{-1}\overset{\text{Lem }\ref{lem:root.conj}}{=}U_{\beta,\,c_\RamiD{0}n+\beta(g)}\overset{\text{Not }\ref{not:root.sub}\eqref{not:root.filtration}}{<}U_{\beta,n}.
\end{equation}
Also $g(T\cap G_{-cn})g^{-1}=T\cap G_{-cn}<T\cap G_{-n}$ \claude{by \Cref{lem:T.norm}}.
\RamiD{Let $U:=\bfU_B(F)$ and $U^-:=\bfU_{B^-}(F)$.}

By \Cref{lem:iva.decomp} for $P=B$
for every $m\in\Z_{>0}$
\begin{equation}\label{lem:T.conj:3}
    \RamiD{G_{-m}=(U^-\cap G_{-m})(T\cap G_{-m})(U\cap G_{-m}).}
    \end{equation}
Applying this with $m=cn$ and conjugating each factor by $g$ gives
\RamiD{
\begin{align*}
gG_{-cn}g^{-1}&=
g(U^-\cap G_{-cn})g^{-1} g(T\cap G_{-cn})g^{-1} g(U\cap G_{-cn})g^{-1}\\
&\overset{\text{Lem }\ref{lem:U.filtration},\,\ref{lem:T.norm}}{<}
\left(\prod_{\beta\in\Phi_{B^-}}gU_{\beta,c_0n}g^{-1}\right)
(T\cap G_{-n})
\left(\prod_{\beta\in\Phi_B}gU_{\beta,c_0n}g^{-1}\right)\\
&\overset{\eqref{lem:T.conj:2}}{<}\Bigl(\prod_{\beta\in\Phi_{B^-}}U_{\beta,n}\Bigr)(T\cap G_{-n})\Bigl(\prod_{\beta\in\Phi_{B}}U_{\beta,n}\Bigr)\\
&\overset{\text{Lem }\ref{lem:U.filtration},\,\eqref{lem:T.conj:3}}{<}G_{-n}.
\end{align*}
}
\end{proof}

\subsection{Effective results in representation theory}\label{ssec:eff.rep}
We will use the following result from \cite{AGKS_eff_rep}.
\begin{theorem}[{\cite[Theorem C]{AGKS_eff_rep}}]\label{thm:mcexp}
        Let $\pi$ be a representation of $G$ of finite length and $\lambda\in \Lambda^+(A_\bfT)$.
        
                Then there are  
                a constant 
                $\mdef{c_{asm}(\pi,\lambda)}\in \Z_{>0}$ and a monic polynomial $$\mdef{\cP_{\pi,\lambda}}=\sum_{i=0}^{\deg\cP_{\pi,\lambda}} a_i x^i\in \C[x]$$ 
                s.t. for any 
                
        \begin{itemize}
            \item $n\in \Z_{>0}$, 
            \item $v\in \pi^{G_{-n}}$,  $w\in \tilde\pi^{G_{-n}}$,
            \item $k>c_{asm}(\pi,\lambda)n$,
        \end{itemize}
we have
$$\sum_{i=0}^{\deg\cP_{\pi,\lambda}}  a_im_{v,w}(\varpi^{(k+i)\lambda})=0.$$        
\end{theorem}

\section{Functions on cones}\label{sec:cones}
Throughout this section we fix a lattice $\mdef{\Lambda}\subset \R^N$ and  a 
 \tdef{cone} $\mdef{L}\subset \Lambda$, i.e. an additive submonoid defined by linear inequalities with integer coefficients.

Let us first motivate the notions and the approach in this section. We were looking for a class of  functions on $L$ that ``satisfy an effective recurrence relation''. We wanted this \RamiB{class} to be closed both w.r.t. products and  w.r.t.  push forward \RamiB{under} proper homomorphism\RamiB{s} of cones.

We did not find a reasonable class like this. However, we did find a class that is closed w.r.t. products -- this is the \RamiB{class given by the} admissibility notion in \Cref{def:adm} below, and a different class that is closed w.r.t. push forward -- these are the functions whose generating function is rational. 

\RamiB{So, we concentrate on the class of admissible functions, and prove that it is closed under pushforward  w.r.t. proper homomorphism of a cone to $\Z_{\geq 0}$.}
\RamiB{In more details}, we do the following: 
\begin{itemize}
    \item Prove that \DimaA{the class of} admissible functions is closed under multiplications (see \Cref{cor:cone.prod} below).
    \item Prove that the generating function of an admissible function is a rational function (see \Cref{lem:height} below).
    \item Observe that the class of functions that have rational generating functions is closed under pushforward  w.r.t. proper homomorphism of cones.
    \item Observe that in the one dimensional case the 2 classes described above coincide.
    \item Deduce that admissible functions are closed under pushforward  w.r.t. proper homomorphism of a cone to $\Z_{\geq 0}$ (see \Cref{cor:pushforward} below).
\end{itemize}
 As usual we prove quantitative versions of all these statements.
 
\begin{definition}\label{def:adm}
$ $
\begin{enumerate}
    \item\label{def:adm:1} For a polynomial $p=\sum_{i=0}^{\deg p} a_ix^i \in \C[x]$ and an integer $n\in \Z_{>0}$, we will call an element $f\in \C^{\Z_{\geq 0}}$
    a \tdef{$(p,n)$-admissible} element if
for any $k>n$ we have $$\sum_{i=0}^{\deg p} a_i f(k+i)=0.$$
    \item A \tdef{recursion datum} on $L$ is a choice of a generating \RamiB{tuple} $\lambda_1,\dots,\lambda_k$ of $L$ together with non-zero polynomials $p_1,\dots,p_k$.
    \item  For a recursion datum $\cR=((\lambda_1,\dots,\lambda_k),(p_1,\dots,p_k))$ and an integer $n\in \Z_{>0}$ we say that $f\in \C^L$ is \tdef{$(\cR,n)$-admissible} if
     for any $\lambda\in L$ and any $i\in\{1,\dots,k\}$ the function $m\mapsto f(\lambda+m\lambda_i)$ is $(p_i,n)$-admissible as an element of  $\C^{\Z_{\geq 0}}$.
\end{enumerate}
    
\end{definition}
\subsection{\Rami{Multiplicativity of admissibility}}\label{ssec:mult.adm}
 \begin{lemma}\label{lem:1d.prod}
     Let  $p_1,p_2 \in \C[x]\RamiD{\smallsetminus 0}$,
     then there are  $p \in \C[x]\RamiD{\smallsetminus 0}$ and $c\in\Z_{>0}$  s.t. for any
     \begin{itemize}
          \item
     $n\in \Z_{>0}$, \item $f_1,f_2\in \C^{\Z_{\geq 0}}$ , s.t. $f_i$ is $(p_i,n)$-admissible (for $i=1,2$),
          \end{itemize}
the function
     $f_1f_2$ is $(p,cn)$-admissible.
 \end{lemma}
  \begin{proof}
 Let $d_i=\deg p_i$, $d=\max(d_1,d_2)$, and write $p_i(x)=x^{e_i}q_i(x)$ with $q_i(0)\neq0$ and $e_i\leq d_i$. 
 \RamiB{Take $c:=d+1$.}
 
 Denote by $S$ the shift \RamiA{operator defined by} $(Sf)(k)=f(k+1)$, and by $\alpha_{i,1},\alpha_{i,2},\dots$ the roots of $q_i$ (i.e.\ the non-zero roots of $p_i$), listed with multiplicity. Put
 $$p(x)=\prod_{a,b}\bigl(x-\alpha_{1,a}\alpha_{2,b}\bigr)^{\,d_1+d_2},$$
 the product over all pairs; $p$ and $c$ depend only on $p_1,p_2$.

 Fix $i\in\{1,2\}$. That $f_i$ is $(p_i,n)$-admissible means $\RamiA{(p_i(S)f_i)}(k)=0$ for $k>n$, i.e.\ $\RamiA{(q_i(S)f_i)}(m)=0$ for $m>n+e_i$. Since $q_i(0)\neq0$, on the half-line $\{m>n+e_i\}$ the sequence $\RamiA{m\mapsto f_i(m)}$ solves a constant-coefficient linear recursion whose characteristic roots are the $\alpha_{i,a}\neq0$; hence
 $$f_i(m)=\sum_a P_{i,a}(m)\,\alpha_{i,a}^{\,m}\qquad(m>n+d),$$
 with \RamiA{polynomials $P_{i,a}$ satisfying} $\deg P_{i,a}<d_i$. 
 Consequently, for $m>n+d$,
 $$f_1(m)f_2(m)=\sum_{a,b}P_{1,a}(m)P_{2,b}(m)\,(\alpha_{1,a}\alpha_{2,b})^{\,m}.$$
 Each summand $R(m)\gamma^m$ with $\deg R<d_1+d_2$ satisfies $(S-\gamma)^{d_1+d_2}\bigl(R(m)\gamma^m\bigr)=0$. So $p(S)(f_1f_2)(k)=0$ for all $k>n+d$. 
 \RamiB{Thus} the function $f_1f_2$ is $(p,cn)$-admissible.
 \end{proof}
\begin{lemma}\label{lem:prop}
     Let  $\cR$ be a recursion datum on $L$ and let $\lambda\in L$.
     Then there is a non-zero $p\in \C[x]$  s.t. for any
     \begin{itemize}
          \item
     $n\in \Z_{>0}$,
     \item $f\in \C^{L}$ , s.t. $f$ is $(\cR,n)$-admissible,
     \item $\lambda_0\in L$,
          \end{itemize}
the function $m\mapsto f(\lambda_0+m\lambda)$ is
     $(p,n)$-admissible as an element of  $\C^{\Z_{\geq 0}}$.
\end{lemma}
\begin{proof}
Write \DimaA{$\cR=((\lambda_1,\dots,\lambda_k),(p_1,\dots,p_k))$ and}  $\lambda=\sum_{i\in I}c_i\lambda_i$ with $I=\{i:c_i>0\}$ and $c_i\in\Z_{>0}$ (the $\lambda_i$ generate $L$). Fix $\lambda_0\in L$ and define $\cF\in\C^{\Z_{\geq0}^I}$ by $\cF\bigl((m_i)_{i\in I}\bigr)=f\bigl(\lambda_0+\sum_{i\in I}m_i\lambda_i\bigr)$, and let $S_i$ be the commuting shift \RamiA{operators defined by} $(S_i\cF)(m)=\cF(m+e_i)$
\RamiB{, where $\{e_i\}_{i\in I}$ is the standard basis of $\Z_{\geq 0}^I$}. For each $i\in I$, fixing the coordinates $m_j$ $(j\neq i)$ the sequence $t\mapsto \cF(m)|_{m_i=t}$ equals $t\mapsto f\bigl((\lambda_0+\sum_{j\neq i}m_j\lambda_j)+t\lambda_i\bigr)$, which is $(p_i,n)$-admissible; that is, $\RamiA{(p_i(S_i)\cF)}(m)=0$ whenever $m_i>n$.

In the finite-dimensional commutative $\C$-algebra $A:=\bigotimes_{i\in I}\C[S_i]/(p_i(S_i))$, of dimension $\prod_{i\in I}\deg p_i$, the element $\prod_{i\in I}S_i^{\,c_i}$ satisfies its monic minimal polynomial $p$; hence, as operators on $\C^{\Z_{\geq0}^I}$,
$$p\Bigl(\textstyle\prod_{i\in I}S_i^{\,c_i}\Bigr)=\sum_{i\in I}r_i(S)\,p_i(S_i),\qquad r_i\in\C[S_j:j\in I],$$
with $p$ (of degree $\leq\prod_i\deg p_i$) depending only on $\cR$ and $\lambda$.

Put $g(m):=f(\lambda_0+m\lambda)=\cF\bigl((c_im)_{i\in I}\bigr)$ and let $S$ be the shift on $\C^{\Z_{\geq0}}$; since $g(m+1)=\bigl(\prod_iS_i^{c_i}\cF\bigr)\bigl((c_im)_i\bigr)$,
$$p(S)g(m)=\Bigl(p\bigl(\textstyle\prod_iS_i^{c_i}\bigr)\cF\Bigr)\bigl((c_im)_i\bigr)=\sum_{i\in I}\bigl(r_i(S)\,p_i(S_i)\cF\bigr)\bigl((c_im)_i\bigr).$$
By the relation above, $p_i(S_i)\cF$ vanishes at every point whose $i$-th coordinate exceeds $n$; as $r_i(S)$ is a combination of shifts $S^b$ ($b\geq0$), which only raise coordinates, $r_i(S)p_i(S_i)\cF$ vanishes there too. Since $c_i\geq1$, the $i$-th coordinate $c_im$ exceeds $n$ once $m>n$; hence $p(S)g(m)=0$ for all $m>n$, i.e.\ $g$ is $(p,n)$-admissible.
\end{proof}
 \begin{cor}\label{cor:cone.prod}
     Let  $\cR_1,\cR_2$   be recursion data on $L$.
     Then there are  a recursion datum $\cR$ on $L$ and $c\in\Z_{>0}$  s.t. for any
     \begin{itemize}
          \item
     $n\in \Z_{>0}$, \item $f_1,f_2\in \C^{L}$ , s.t. $f_i$ is $(\cR_i,n)$-admissible (for $i=1,2$),
          \end{itemize}
the function
     $f_1f_2$ is $(\cR,cn)$-admissible.
 \end{cor}
  \begin{proof}
 Let $\lambda_1,\dots,\lambda_k$ be the generating set of $\cR_1$, with polynomials $p_1,\dots,p_k$. Fix $i$. Along $\lambda_i$, $f_1$ is $(p_i,n)$-admissible by definition, and by \Cref{lem:prop} applied to $\cR_2$ and the vector $\lambda_i$ there is a polynomial $p_i'$ (depending only on $\cR_2,\lambda_i$) such that $m\mapsto f_2(\lambda_0+m\lambda_i)$ is $(p_i',n)$-admissible for every $\lambda_0\in L$. By \Cref{lem:1d.prod} there are a polynomial $r_i$ and $c_i\in\Z_{>0}$ (depending only on $p_i,p_i'$) such that the product $m\mapsto(f_1f_2)(\lambda_0+m\lambda_i)=f_1(\lambda_0+m\lambda_i)\,f_2(\lambda_0+m\lambda_i)$ is $(r_i,c_in)$-admissible.

 Take $\cR$ to be the datum with generating set $\lambda_1,\dots,\lambda_k$ and polynomials $r_1,\dots,r_k$, and $c:=\max_i c_i$. For every $\lambda_0\in L$ and every $i$ the sequence $m\mapsto(f_1f_2)(\lambda_0+m\lambda_i)$ is $(r_i,c_in)$-admissible, hence $(r_i,cn)$-admissible; thus $f_1f_2$ is $(\cR,cn)$-admissible.
 \end{proof}
\subsection{\Rami{Extension of admissibility from a subcone of finite index}}\label{ssec:ext.adm}
\begin{lemma}\label{lem:fin.index}
    Let $\Lambda_0<\Lambda$ be a subgroup of finite index and let $L_0:=L\cap \Lambda_0$.     
    Then for any recursion datum $\cR_0$ on $L_0$ and any $Z\subset L$ satisfying 
    $L_0+Z=L$ there are a recursion datum $\cR$ on $L$ and $c\in \Z_{>0}$ s.t. for any
    \begin{itemize}
        \item $n\in \Z_{>0}$,
        \item $f\in \C^L$ s.t. for any $z\in Z$ the function $\lambda\mapsto f(z+\lambda)$ is $(\cR_0,n)$-admissible as an element of $\C^{L_0}$,
    \end{itemize}
    the function $f$ is $(\cR,cn)$-admissible.
\end{lemma}
\begin{proof}
Put $N_0:=[\Lambda:\Lambda_0]$ and fix a generating set $\lambda_1,\dots,\lambda_k$ of $L$. Since $\Lambda/\Lambda_0$ is a group of order $N_0$ we have $N_0\Lambda\subset \Lambda_0$, so, $L$ being a monoid, $N_0\lambda_j\in L\cap \Lambda_0=L_0$ for every $j$. Note that $\Lambda_0$ is a lattice and that $L_0$ is a cone in it, being cut out in $\Lambda_0$ by the linear inequalities that cut out $L$ in $\Lambda$; so the \RamiB{previous} results of this section apply to the pair $(\Lambda_0,L_0)$ as well.

For every $j$, \Cref{lem:prop}, applied to the cone $L_0$, the recursion datum $\cR_0$ and the element $N_0\lambda_j\in L_0$, provides a non-zero polynomial $q_j=\sum_{l=0}^{d_j}a^{(j)}_lx^l$ s.t. for any $n'\in \Z_{>0}$, any $(\cR_0,n')$-admissible $h\in \C^{L_0}$ and any $\mu\in L_0$ the function $t\mapsto h(\mu+tN_0\lambda_j)$ is $(q_j,n')$-admissible. Put
$$p_j(x):=q_j(x^{N_0}),\qquad \cR:=\bigl((\lambda_1,\dots,\lambda_k),(p_1,\dots,p_k)\bigr),\qquad c:=2N_0 .$$

Let $n$ and $f$ be as in the formulation, let $\lambda_0\in L$ and $j\in\{1,\dots,k\}$, and put $F(m):=f(\lambda_0+m\lambda_j)$. Let $m>cn$ and write $m=r+N_0t$ with $r\in \{0,\dots,N_0-1\}$ and $t\in \Z_{\geq 0}$; then $N_0t=m-r>2N_0n-N_0\geq N_0n$, so $t>n$. As $L_0+Z=L$ there are $z\in Z$ and $\mu\in L_0$ with $\lambda_0+r\lambda_j=z+\mu$. The function $h\in \C^{L_0}$ given by $h(\lambda):=f(z+\lambda)$ is $(\cR_0,n)$-admissible, hence $t'\mapsto h(\mu+t'N_0\lambda_j)$ is $(q_j,n)$-admissible and therefore
$$\sum_{l=0}^{d_j}a^{(j)}_lF(m+N_0l)=\sum_{l=0}^{d_j}a^{(j)}_lh\bigl(\mu+(t+l)N_0\lambda_j\bigr)\overset{\ref{def:adm}\eqref{def:adm:1}}{=}0 ,$$
the first equality because $\lambda_0+(m+N_0l)\lambda_j=z+\mu+(t+l)N_0\lambda_j$, and the second one being the recursion of $q_j$ at the point $t>n$. Since $p_j(x)=\sum_{l}a^{(j)}_lx^{N_0l}$, this says that $F$ is $(p_j,cn)$-admissible. As $\lambda_0$ and $j$ were arbitrary, $f$ is $(\cR,cn)$-admissible.
\end{proof}
\RamiE{
The following lemmas are straightforward:
\begin{lemma}\claude{\label{lem:many.data}}
    For any \claude{finite} collection of \claude{recursion} data $\{\cR_i\}_{i\in I}$ on $L$ there is a \claude{recursion} datum $\cR$ on $L$ and a constant $c$ s.t. if  $f\in \C^L$ is $(\cR_i,n)$-\claude{admissible} for some $i\in I$ \claude{and $n\in\Z_{>0}$} then it is $(\cR,cn)$-\claude{admissible}.
\end{lemma}
}
\RamiE{
\begin{lemma}\claude{\label{lem:scale.data}}
        For any {recursion} datum $\cR$ on $L$  and integer $e\in \Z_{>0}$ there is a {recursion} datum $\mdef{\cR^e}$ on $eL$ \claude{s.t.} for any $n\in\Z_{>0}$, if  $f\in \C^L$ is $(\cR,n)$-{admissible} then $f|_{eL}$ is $(\cR^e,n)$-{admissible}.

\end{lemma}
}

\subsection{\Rami{Pushforward of admissible functions}}\label{ssec:push.adm}

 \begin{notation}
 $ $
 \begin{itemize}
     \item For an element $\lambda\in L$ we denote by $\mdef{\delta_\lambda}$ the corresponding element in the monoid algebra $\mdef{\C[L]}$.
     \item 
     Denote by $\mdef{\cI_L}\subset \C[L]$ the ideal generated by non-invertible elements in $L$.

     \item Denote the completion $\widehat{\C[L]}_{\cI_L}$ by $\mdef{\C[[L]]}$.
         \end{itemize} 
 \end{notation}

The following \DimaA{lemma} is standard.
\begin{lemma}
Assume that the only invertible element in $L$ is $0$. Then
$ $
    \begin{itemize}
        \item $\cI_L$ is a maximal ideal.
        \item $\C^L$ is canonically isomorphic \RamiA{(as a vector space)} to $\C[[L]]$.
    \end{itemize}
\end{lemma}
\begin{notation}
    Assume that the only invertible element in $L$ is $0$. Then for $f\in \C^L$ we denote the corresponding element in $\C[[L]]$ by $\mdef{\hat{f}}$. 

    We also denote the inverse map by $\mdef{\cC}:\C[[L]]\to \C^L$.
\end{notation}

\begin{notation}
$ $
\begin{itemize}
    \item 
    For a polynomial $p=\sum_{i=0}^{d}a_ix^i$ with $d=\deg p$ we write $\mdef{\tilde p}(x)=x^{d}p(1/x)=\sum_{i=0}^{d}a_{d-i}x^i$. 
    \item We fix integer linear forms $\mdef[\nu]{\nu_1,\dots,\nu_t}$ on $\Lambda$ with
    $$L=\{x\in\Lambda:\langle\nu_a,x\rangle\geq0\text{ for all }a\},$$
    which exist by the definition of a cone. \item
    For $\lambda\in L$ and $z\in L$ we set
    $$\mdef{\mathrm{ord}_\lambda(z)}=\max\{m\in\Z_{\geq0}:z-m\lambda\in L\}\in\Z_{\geq0}\cup\{\infty\}.$$
\end{itemize}
\end{notation}

\begin{lemma}\label{lem:one-gen}
\RamiE{Assume that the only invertible element in $L$ is $0$.}
Let 
\begin{itemize}    
    \item 
    $\cR=((\lambda_1,\dots,\lambda_k),(p_1,\dots,p_k))$ be a recursion datum,
    \item \RamiB{$i\in \{1,\dots, k\}$,}
    \item $D_i=\deg p_i$, \item $n\in\Z_{>0}$,
    \item $f\in\C^L$ be an $(\cR,n)$-admissible function. 
    \end{itemize} 
    Then $u:=\cC\bigl(\hat f\,\tilde p_i(\delta_{\lambda_i})\bigr)$ satisfies $u(z)=0$ for every $z\in L$ with $\mathrm{ord}_{\lambda_i}(z)>n+D_i$.
\end{lemma}
\begin{proof}
    Write $p_i=\sum_{r=0}^{D_i}a_r x^r$, so
    \begin{align*}
    \tilde p_i(\delta_{\lambda_i})=\sum_{s=0}^{D_i}a_{D_i-s}\,\delta_{s\lambda_i},\qquad
    u(z)=\sum_{s=0}^{D_i}a_{D_i-s}f(z-s\lambda_i)
    \end{align*}
    (with $f$ extended by $0$ off $L$). Since $\mathrm{ord}_{\lambda_i}(z)>n+D_i$, choose an integer $M>n+D_i$ with $y:=z-M\lambda_i\in L$. For $0\leq s\leq D_i$ then $z-s\lambda_i=y+(M-s)\lambda_i\in L$ (as $M-s\geq0$ and $\lambda_i\in L$), and since $t\mapsto f(y+t\lambda_i)$ is $(p_i,n)$-admissible,
    \begin{align*}
    u(z)&=\sum_{s=0}^{D_i}a_{D_i-s}\,f(y+(M-s)\lambda_i)\\
    &=\sum_{r=0}^{D_i}a_r\,f\bigl(y+((M-D_i)+r)\lambda_i\bigr)
    \overset{\ref{def:adm}\eqref{def:adm:1}}{=}0,
    \end{align*}
    the last equality being the defining recurrence at index $M-D_i>n$.
\end{proof}
\begin{lemma}\label{lem:reentry}
    Let $\RamiB{\lambda,\mu\in L}$ with $\langle\nu_a,\RamiB{\lambda}\rangle>0$ for some $a\RamiB{\in \{1,\dots, t\}}$. Then there is $\RamiB{C_{\lambda\mu}}\in\Z_{\geq0}$ such that $\mathrm{ord}_{\RamiB{\lambda}}(z)\leq\mathrm{ord}_{\RamiB{\lambda}}(z-\RamiB{\mu})+\RamiB{C_{\lambda\mu}}$ for all $z\in L$ with $z-\RamiB{\mu}\in L$.
\end{lemma}
\begin{proof}
    For $z\in L$ and $m\in\Z_{\geq0}$ one has $z-m\RamiB{\lambda}\in L$ iff $\langle\nu_a,z\rangle\geq m\langle\nu_a,\RamiB{\lambda}\rangle$ for all $a\RamiB{\in \{1,\dots,t\}}$. The constraints with $\langle\nu_a,\RamiB{\lambda}\rangle=0$ hold automatically, as $\langle\nu_a,z\rangle\geq0$. Hence
    $$\mathrm{ord}_{\RamiB{\lambda}}(z)=\min_{a:\langle\nu_a,\RamiB{\lambda}\rangle>0}\Big\lfloor\frac{\langle\nu_a,z\rangle}{\langle\nu_a,\RamiB{\lambda}\rangle}\Big\rfloor.$$
    Set $\RamiB{C_{\lambda\mu}}=\max_{a:\langle\nu_a,\RamiB{\lambda}\rangle>0}\big\lceil\langle\nu_a,\RamiB{\mu}\rangle/\langle\nu_a,\RamiB{\lambda}\rangle\big\rceil\in\Z_{\geq0}$. Since $$\langle\nu_a,z\rangle=\langle\nu_a,z-\RamiB{\mu}\rangle+\langle\nu_a,\RamiB{\mu}\rangle,$$ for each $a$ with $\langle\nu_a,\RamiB{\lambda}\rangle>0$ we have $$\lfloor\langle\nu_a,z\rangle/\langle\nu_a,\RamiB{\lambda}\rangle\rfloor\leq\lfloor\langle\nu_a,z-\RamiB{\mu}\rangle/\langle\nu_a,\RamiB{\lambda}\rangle\rfloor+\RamiB{C_{\lambda\mu}}.$$ Taking the minimum over such $a$ gives the claim.
\end{proof}

\begin{lemma}\label{lem:height}
    Assume the only invertible element in $L$ is $0$. For any recursion datum $\cR=((\lambda_1,\dots,\lambda_k),(p_1,\dots,p_k))$ there is a monoid homomorphism $h:L\to(\R_{\geq0},+)$ with $Ker(h)=0$ such that for any
\begin{itemize}
    \item $n\in \Z_{>0}$
    \item an $(\cR,n)$-admissible $f\in \C^L$
\end{itemize}
 we have
    $$\supp\left(\cC\left(\hat{f}\prod_{i=1}^k \tilde p_i(\delta_{\lambda_i})\right)\right)\subset h^{-1}([0,n]).$$
\end{lemma}
 \begin{proof}
 We may assume each $\lambda_i\neq0$. As $0$ is the only invertible element of $L$, no nonzero $x\in L$ satisfies $\langle\nu_a,x\rangle=0$ for all $a$ (else $-x\in L$, so $x$ is a unit); hence for every $i$ there is $a$ with $\langle\nu_a,\lambda_i\rangle>0$, and \Cref{lem:reentry} provides constants $\RamiB{C_{ij}:=C_{\lambda_i\lambda_j}}$. Set $D_i=\deg p_i$, $K_i=1+D_i+\sum_{j}D_jC_{ij}$, $N_i=n+D_i+\sum_{j\neq i}D_jC_{ij}$ (so $N_i\leq nK_i$ for $n\geq1$), and $g=\cC\bigl(\hat f\prod_i\tilde p_i(\delta_{\lambda_i})\bigr)$.

 Fix $i$ and factor $\prod_l\tilde p_l(\delta_{\lambda_l})=\tilde p_i(\delta_{\lambda_i})\,Q_i$ with $Q_i=\prod_{j\neq i}\tilde p_j(\delta_{\lambda_j})\in\C[L]$ supported on $B_i=\{\sum_{j\neq i}l_j\lambda_j:0\leq l_j\leq D_j\}$. Writing $u_i=\cC\bigl(\hat f\,\tilde p_i(\delta_{\lambda_i})\bigr)$ we get $g(z)=\sum_{\mu\in B_i}(Q_i)_\mu\,u_i(z-\mu)$. If $g(z)\neq0$, some term is nonzero: there is $\mu=\sum_{j\neq i}l_j\lambda_j\in B_i$ with $z-\mu\in L$ and $u_i(z-\mu)\neq0$, so $\mathrm{ord}_{\lambda_i}(z-\mu)\leq n+D_i$ by \Cref{lem:one-gen}. Applying \Cref{lem:reentry} $l_j$ times for each $j\neq i$ (the intermediate points lie in $L$),
 \begin{align*}
 \mathrm{ord}_{\lambda_i}(z)&\overset{\text{Lem }\ref{lem:reentry}}{\leq}\mathrm{ord}_{\lambda_i}(z-\mu)+\sum_{j\neq i}l_jC_{ij}\\
 &\overset{\text{Lem }\ref{lem:one-gen}}{\leq} n+D_i+\sum_{j\neq i}D_jC_{ij}=N_i.
 \end{align*}
 As $i$ was arbitrary, $\supp(g)\subset\{z\in L:\mathrm{ord}_{\lambda_i}(z)\leq N_i\text{ for all }i\}$. For $z=\sum_r c_r\lambda_r\in L$ ($c_r\in\Z_{\geq0}$) one has $z-c_r\lambda_r=\sum_{s\neq r}c_s\lambda_s\in L$, so $\mathrm{ord}_{\lambda_r}(z)\geq c_r$; thus every $z\in\supp(g)$ has $c_r\leq N_r$ for all $r$ in any such representation, and $\supp(g)$ lies in the finite set $\{\sum_r c_r\lambda_r:0\leq c_r\leq N_r\}$.

 As $0$ is the only invertible element of $L$, there is a linear functional $\xi$ on $\Lambda$ with $\xi>0$ on $L\setminus\{0\}$. Put $A=\sum_r K_r\xi(\lambda_r)>0$ and $h=\xi/A$, a monoid homomorphism with $Ker(h)=0$. For $z=\sum_r c_r\lambda_r\in\supp(g)$ (so $c_r\leq N_r\leq nK_r$),
 $$h(z)=\frac1A\sum_r c_r\xi(\lambda_r)\leq\frac1A\sum_r nK_r\xi(\lambda_r)=n.$$
 Hence $\supp(g)\subset h^{-1}([0,n])$.
 \end{proof}

 
 \begin{cor}\label{cor:pushforward}
 Let $\phi:L\to \Z_{\geq 0}$ be a monoid homomorphism with trivial kernel. Then for any recursion datum $\cR=((\lambda_1,\dots,\lambda_k),(p_1,\dots,p_k))$ there are a nonzero polynomial $p$ and $c\in\Z_{> 0}$  such that for any $n\in \Z_{>0}$ and any $(\cR,n)$-admissible $f\in \C^L$ the function  $\phi_*(f)$ is $(p,cn)$-admissible.
  \end{cor}
 \begin{proof}
 \RamiE{We may assume each $\lambda_i\neq0$.}
 Since $Ker(\phi)=0$ the only unit of $L$ is $0$, and $w_i:=\phi(\lambda_i)\in\Z_{>0}$. Keep the notation $D_i,C_{ij},N_i,K_i$ of \Cref{lem:height}. By its proof, $g:=\cC(\hat fP)$ with $P=\prod_i\tilde p_i(\delta_{\lambda_i})$ is supported in the finite set $\{\sum_r c_r\lambda_r:0\leq c_r\leq N_r\}$, and there
 $$\phi(z)=\sum_r c_rw_r\leq\sum_r N_rw_r\leq n\sum_r K_rw_r=cn,\qquad c:=\sum_r K_rw_r\in\Z_{>0},$$
 using $N_r\leq nK_r$ ($n\geq1$).

 As $Ker(\phi)=0$ and $L$ is a strict cone, $\phi$ has finite fibres, so $\phi_*:\C[[L]]\to\C[[\Z_{\geq0}]]$ is a well-defined $\C$-algebra homomorphism with $\phi_*\hat f=\widehat{\phi_*f}$. Since $\phi_*(\delta_{\lambda_i})=\delta^{w_i}$ (where $\delta=\delta_1$ generates $\C[\Z_{\geq0}]$), applying $\phi_*$ to $\hat g=\hat fP$ gives
 $$\widehat{\phi_*f}\cdot q(\delta)=\phi_*\hat g,\qquad q(x):=\prod_{i=1}^k\tilde p_i(x^{w_i}).$$
 The right-hand side is the pushforward of the finitely supported $g$, hence vanishes above $\max_{\supp g}\phi\leq cn$. Writing $\cF=\phi_*f$ and $q(x)=\sum_{l=0}^{D_q}q_lx^l$ with $D_q:=\sum_iw_iD_i\geq\deg q$, the identity reads $\sum_{l=0}^{D_q} q_l\cF(m-l)=0$ for all $m>cn$. Substituting $l=D_q-j$ and $k=m-D_q$ turns this into $\sum_{j=0}^{D_q} q_{D_q-j}\cF(k+j)=0$ for all $k>cn-D_q$, in particular for all $k>cn$. Thus $\phi_*f$ is $(p,cn)$-admissible for the nonzero polynomial $p(x)=x^{D_q} q(1/x)$.
 \end{proof}
 \begin{cor}\label{cor:acc_pushforward}
 Let $\phi:L\to \Z_{\geq 0}$ be a monoid homomorphism with trivial kernel. Then for any recursion datum $\cR$ there are a nonzero polynomial $p$ and $c\in\Z_{> 0}$  such that for any $n\in \Z_{>0}$ and any $(\cR,n)$-admissible $f\in \C^L$ the function   
 $$k\mapsto \sum_{x\in L; \phi(x)\leq k} f(x)$$ is $(p,cn)$-admissible.
  \end{cor}
 \section{Admissibility of matrix coefficients}\label{sec:adm.mat}
Fix an irreducible representation $\pi$ of $G$. 
Let $\mdef{\chi}$ be the central character of $\pi$.

In this section we prove 
that for a matrix coefficient $m$ of $\pi$ and a test function  $f$ the function $\cS:\Lambda^+(A_\bfT)\to \C$ given by $$\cS(\lambda)=\langle Ad(\varpi^\lambda)(m),f\rangle$$
is  $(\cR,n)$-admissible for some fixed $\cR$ (depending on $\pi$) and for $n$ that depends linearly on the level of $m$ and $f$. 
See \Cref{cor:S.adm}.
This is an application of the results of \cite{AGKS_eff_rep} and it is used in the proof of \Cref{thm:recur}.

\begin{notation}
    Denote by $\mdef{\Lambda^0(A_{\bfT})}$  the set of invertible elements in 
    $\Lambda^+(A_{\bfT})$. \RamiE{Let        $\mdef{\Lambda_{ad}(A_{\bfT})}=
    \Lambda(A_{\bfT})/\Lambda^0(A_{\bfT})$
     and} 
    $\mdef{\Lambda^+_{ad}(A_{\bfT})}=
    \Lambda^+(A_{\bfT})/\Lambda^0(A_{\bfT})$.
\end{notation}
\Rami{
\begin{lemma}\label{lem:gen.cone.fac}
    Let $\Lambda\subset \R^N$ and $L\subset \Lambda$ be as in \Cref{sec:cones}. Assume that $L$ is given by the inequalities $f_1\geq 0,\dots,f_k\geq 0$, where $f_1,\dots,f_k:\Lambda\to \Q$ are homomorphisms, i.e.
    $$L=\{\lambda\in \Lambda\,|\,\forall i\in \{1,\dots,k\} \text{ we have } f_i(\lambda)\geq 0\}.$$
    Let
    $$\Lambda^0:=\{\lambda\in \Lambda\,|\,\forall i\in \{1,\dots,k\} \text{ we have } f_i(\lambda)= 0\}.$$
    Then $\Lambda^0$ is the group of invertible elements of the monoid $L$, the group $\Lambda/\Lambda^0$ is a lattice, and the image $L/\Lambda^0$ of $L$ in $\Lambda/\Lambda^0$ is a cone whose only invertible element is $0$.
\end{lemma}
\begin{proof}
WLOG we may assume that the images of $f_i$ are  in $\Z$.
    Let $\varphi:=(f_1,\dots,f_k):\Lambda\to \Z^k$. Then $\Lambda^0=\ker \varphi$ and $L=\varphi^{-1}(\Z_{\geq 0}^k)$. In particular $\Lambda^0\subset L$ and $L$ is a union of cosets of $\Lambda^0$, so $L/\Lambda^0$ is a submonoid of $\Lambda/\Lambda^0$ and $L=\varphi^{-1}(\varphi(L))$.

    An element $\lambda\in L$ is invertible in $L$ if and only if $-\lambda\in L$, i.e. if and only if $f_i(\lambda)\geq 0$ and $f_i(\lambda)\leq 0$ for every $i$, i.e. if and only if $\lambda\in \Lambda^0$.

    The map $\varphi$ induces an injective homomorphism $\bar\varphi:\Lambda/\Lambda^0\to \Z^k$, so $\Lambda/\Lambda^0$ is isomorphic to its image, which is a lattice in $\R^k$. Under this isomorphism
    $$\bar\varphi(L/\Lambda^0)=\varphi(L)=\bar\varphi(\Lambda/\Lambda^0)\cap \Z^k_{\geq 0},$$
    where the second equality holds since $L=\varphi^{-1}(\varphi(L))$. Hence $L/\Lambda^0$ is defined inside the lattice $\Lambda/\Lambda^0$ by the inequalities $x_1\geq 0,\dots,x_k\geq 0$, where $x_1,\dots,x_k$ are the coordinates of $\R^k$; so it is a cone. Finally, if $x\in L/\Lambda^0$ and $-x\in L/\Lambda^0$ then every coordinate of $\bar\varphi(x)$ vanishes, so $x=0$.
\end{proof}
}

\begin{lemma}\label{lem:Lam.ad.cone}
    We have
    $$\Lambda^0(A_{\bfT})=\{\nu\in \Lambda(A_{\bfT})\,|\,\forall \beta\in \Phi_{\Rami{\bfB}} \text{ we have } \beta(\nu)=0\}.$$
    The group $\Lambda(A_{\bfT})/\Lambda^0(A_{\bfT})$ is a lattice, and $\Lambda^+_{ad}(A_{\bfT})$ is a cone in it whose only invertible element is $0$ \claude{and which generates $\Lambda(A_{\bfT})/\Lambda^0(A_{\bfT})$}.
\end{lemma}
\begin{proof}
    By \Cref{not:settings}\eqref{it:M+},\eqref{it:Lam+},\eqref{it:Lam+A} we have
        \Rami{$$\Lambda^+(A_{\bfT})=\{\nu\in \Lambda(A_{\bfT})\,|\,\forall \beta\in \Phi_{{\bfB}} \text{ we have } \beta(\nu)\geq 0\},$$}
    and $\Phi=\Phi_\bfB\sqcup(-\Phi_\bfB)$ by \Cref{not:settings}\eqref{it:Phi.P}. Hence $\nu$ is invertible in the monoid $\Lambda^+(A_{\bfT})$ iff both $\nu$ and $-\nu$ lie in $\Lambda^+(A_{\bfT})$, i.e.\ iff $\beta(\nu)=0$ for every $\beta\in\Phi$. This proves the first assertion. 
    \Rami{The rest of the Lemma follows now from \Cref{lem:gen.cone.fac}}    
\end{proof}

\begin{lemma}\label{lem:Lam0.center}$\varpi^{\Lambda^0(A_{\bfT})}=\varpi^{\Lambda^+(A_{\bfT})}\cap \cZ(G)$.
\end{lemma}
\begin{proof}
    By \Cref{not:settings}\eqref{it:M+}, \eqref{it:Lam+} and \eqref{it:Lam+A} we have         \Rami{$$\Lambda^+(A_{\bfT})=\{\nu\in \Lambda(A_{\bfT})\,|\,\forall \beta\in \Phi_{{\bfB}} \text{ we have } \beta(\nu)\geq 0\},$$}
    and since every element of $\Phi_\bfB$ is a non-negative integral combination of the simple roots $\Delta$ of \Cref{not:settings}\eqref{it:Delta},
    $$\Lambda^+(A_{\bfT})=\{\lambda\in\Lambda(A_{\bfT}):\alpha(\lambda)\geq0\ \forall\alpha\in\Delta\}.$$
    For $\lambda\in\Lambda(A_{\bfT})$ the element $\varpi^\lambda=\cE(\lambda)\in A_\bfT$ centralizes $T$ and normalizes each relative root group $U_\alpha$, acting on the weight spaces of $\alpha$ and of $2\alpha$ in $\Lie(U_\alpha)$ (\Cref{not:root.sub}\eqref{not:root.grp}) by the scalars $\varpi^{\alpha(\lambda)}$ and $\varpi^{2\alpha(\lambda)}$ respectively; both are $1$ iff $\alpha(\lambda)=0$, and then $\varpi^\lambda$ centralizes $U_\alpha$. As the $U_\alpha$ and $T$ generate $G$,
    \begin{align*}
    \varpi^\lambda\in \cZ(G)&\iff \alpha(\lambda)=0\ \text{for every relative root }\alpha\\
    &\iff \alpha(\lambda)=0\ \forall\alpha\in\Delta.
    \end{align*}
    On the other hand, \Rami{by \Cref{lem:Lam.ad.cone},} $\lambda\in\Lambda^0(A_{\bfT})$, iff $\alpha(\lambda)=0\ \forall\alpha\in\Delta$. Thus, for $\lambda\in\Lambda^+(A_{\bfT})$, $\varpi^\lambda\in \cZ(G)$ iff $\lambda\in\Lambda^0(A_{\bfT})$, which is the asserted equality of subsets of $G$.
\end{proof}

So we can pull $\chi$ to $\Lambda^0(A_{\bfT})$. We will  denote the pullback by $\chi$ as well. Fix a \RamiC{(homomorphic)} extension of it to $\Lambda^+(A_{\bfT})$ and denote it by $\mdef{\tilde\chi}$.

\begin{notation}   
For $\lambda\in \Lambda^+(A_{\bfT})$ and $v\in \pi, w\in \tilde\pi$ and $f\in \cH(G)$, denote
\begin{itemize}
\item $\mdef{\bar m_{v,w}}(\lambda)=\tilde\chi^{-1}(\lambda) m^\pi_{v,w}(\varpi^\lambda)$.
    \item 
     $\mdef{\cS_{v,w,f}}(\lambda)=\langle Ad(\varpi^\lambda)( m_{ v,w}),f\rangle=\langle m_{\varpi^\lambda v,\varpi^\lambda w},f\rangle$.
\end{itemize}
We consider both $\bar m_{v,w}$ and $\cS_{v,w,f}$ as elements in 
$\C^{\Lambda^+_{ad}({A_{\bfT}})}$.
\end{notation}

\begin{lemma}[{\cite[Lemma V.5.3.]{Ren}}\footnote{\Cref{lem:iva.decomp}, \RamiD{\Cref{lem:T.norm},  \Cref{lem:U.filtration}\eqref{lem:U.filtration:3} and \Cref{lem:parahoric}\eqref{lem:parahoric:norm} ensure} that the conditions of \cite[Lemma V.5.3.]{Ren} are satisfied.}]\label{lem:tau.semigroup}
Fix $n\in\Z_{>0}$ and for $\lambda\in\Lambda^+(A_\bfT)$ set
$\mdef{\tau_\lambda}:=e_{G_{-n}}\delta_{\varpi^\lambda}e_{G_{-n}}\in\cH_{G_{-n}}(G)$.
\DimaA{Then }for all $\lambda,\mu\in\Lambda^+(A_\bfT)$ we have
$\tau_\lambda \tau_\mu=\tau_{\lambda+\mu}$.
\end{lemma}

\Cref{thm:mcexp} gives us:
\begin{cor}\label{cor:mbar.adm}
    There is a recursion datum $\cR$ on $\Lambda^+_{ad}({A_{\bfT}})$ and a constant $\mdef{c_{asm}}$ s.t. for any $n\in \Z_{>0}$ and any $v\in \pi^{G_{-n}}$  and  $w\in \tilde\pi^{G_{-n}}$, the element  $\bar m_{v,w}$ is $(\cR,c_{asm}n)$-admissible.
\end{cor}
\begin{proof}
    Fix representatives $\lambda_1,\dots,\lambda_r\in\Lambda^+(A_{\bfT})$ of a finite generating set of the monoid $\Lambda^+_{ad}(A_{\bfT})$. For each $i$ let $\cP_{\pi,\lambda_i}=\sum_{j=0}^{d_i}a^{(i)}_jx^j$ (monic, $d_i=\deg\cP_{\pi,\lambda_i}$) and $c_{asm}(\pi,\lambda_i)$ be as in \Cref{thm:mcexp}; put $\zeta_i:=\tilde\chi^{-1}(\lambda_i)\in\C^\times$, $p_i(x):=\sum_{j=0}^{d_i}a^{(i)}_j\zeta_i^{-j}x^j$ (nonzero, as $a^{(i)}_{d_i}=1$ and $\zeta_i\neq0$), $$c_{asm}:=\max_ic_{asm}(\pi,\lambda_i) \text{ and } \cR:=((\lambda_1,\dots,\lambda_r),(p_1,\dots,p_r)).$$
    Fix $n\in\Z_{>0}$, $v\in\pi^{G_{-n}}$, $w\in\tilde\pi^{G_{-n}}$, a base $\mu\in\Lambda^+(A_{\bfT})$ and an index $i$. As $e_{G_{-n}}$ fixes $v$ and $w$, for every $\nu\in\Lambda^+(A_{\bfT})$ we have $m_{v,w}(\varpi^\nu)=\langle\pi(\tau_\nu)v,w\rangle$, with $\tau_\nu$ as in \Cref{lem:tau.semigroup}. Set $T_i:=\pi(\tau_{\lambda_i})\in\End(\pi^{G_{-n}})$ and $v':=\pi(\tau_\mu)v\in\pi^{G_{-n}}$. By \Cref{lem:tau.semigroup}, $\tau_{\mu+k\lambda_i}=\tau_\mu \cdot \tau_{\lambda_i}^{k}$, so $\pi(\tau_{\mu+k\lambda_i})=T_i^k\,\pi(\tau_\mu)$ and
    $$m_{v,w}(\varpi^{\mu+k\lambda_i})\overset{\text{Lem }\ref{lem:tau.semigroup}}{=}\langle T_i^k v',w\rangle=m_{v',w}(\varpi^{k\lambda_i}).$$
    Since $v'\in\pi^{G_{-n}}$ and $w\in\tilde\pi^{G_{-n}}$ (both at level $n$, \emph{independently of $\mu$}), \Cref{thm:mcexp} applied to $v',w$ and the cocharacter $\lambda_i$ shows that $k\mapsto m_{v',w}(\varpi^{k\lambda_i})$ is $(\cP_{\pi,\lambda_i},c_{asm}(\pi,\lambda_i)n)$-admissible.

    Now $\bar m_{v,w}(\mu+k\lambda_i)=\tilde\chi^{-1}(\mu+k\lambda_i)\,m_{v,w}(\varpi^{\mu+k\lambda_i})=\tilde\chi^{-1}(\mu)\,\zeta_i^{\,k}\,m_{v',w}(\varpi^{k\lambda_i})$. In general, if $b$ is $(\cP,N)$-admissible with $\cP=\sum_ja_jx^j$ and $\zeta\in\C^\times$, then $\sum_j(a_j\zeta^{-j})\,\zeta^{k+j}b(k+j)=\zeta^k\sum_ja_jb(k+j)=0$ for $k>N$, so $k\mapsto\zeta^kb(k)$ is $(\sum_ja_j\zeta^{-j}x^j,N)$-admissible; scaling by the constant $\tilde\chi^{-1}(\mu)$ preserves this. Hence $k\mapsto\bar m_{v,w}(\mu+k\lambda_i)$ is $(p_i,c_{asm}(\pi,\lambda_i)n)$-admissible, thus $(p_i,c_{asm}n)$-admissible as $c_{asm}\geq c_{asm}(\pi,\lambda_i)$. As $\mu$ and $i$ were arbitrary, $\bar m_{v,w}$ is $(\cR,c_{asm}n)$-admissible.
\end{proof}

\begin{lemma}
    For $w'\in\tilde\pi$ and $v'\in\pi$ set $\mdef{\check m_{w',v'}}(\lambda):=\tilde\chi(\lambda)\,m^{\tilde\pi}_{w',v'}(\varpi^\lambda)$, the twisted matrix coefficient of the contragredient $\tilde\pi$ (like $\bar m$ it descends to an element of $\C^{\Lambda^+_{ad}(A_{\bfT})}$, since $\varpi^{\Lambda^0(A_{\bfT})}\subset \cZ(G)$ and $\tilde\pi$ has central character $\chi^{-1}$). Then for any $n\in \Z_{>0}$ and any $v\in \pi^{G_{-n}}$, $w\in \tilde\pi^{G_{-n}}$ and $f\in \cH_{G_{-n}}(G)$ there are
    \begin{itemize}
        \item $w_1,\dots, w_k\in \tilde \pi^{G_{-n}}$
        \item $v_1,\dots, v_k\in \pi^{G_{-n}}$
        \item $a_{pq}\in \C$ for $p,q\in \{1,\dots, k\}$
    \end{itemize}
    s.t.     $$\cS_{v,w,f}=\sum_{p,q=1}^k  a_{pq}\,  \bar m_{v,w_p} \cdot \check m_{w,v_q}.$$
\end{lemma}
\begin{proof}
    As $\pi$ is admissible, $\pi^{G_{-n}}$ is finite dimensional; fix a basis $e_1,\dots,e_k$ of it and the dual basis $e_1^\vee,\dots,e_k^\vee\in\tilde\pi^{G_{-n}}$ (so $\langle e_p,e_q^\vee\rangle=\delta_{p,q}$), and set $w_p:=e_p^\vee$, $v_q:=e_q$ and $a_{pq}:=\langle\pi(f)e_p,e_q^\vee\rangle$.

    By definition $\cS_{v,w,f}(\lambda)=\langle m_{\varpi^\lambda v,\varpi^\lambda w},f\rangle=\langle\pi(f)\pi(\varpi^\lambda)v,\tilde\pi(\varpi^\lambda)w\rangle$. Writing $e:=e_{G_{-n}}$ for the idempotent projecting onto level-$n$ vectors, $f\in\cH_{G_{-n}}(G)$ gives $\pi(f)=\pi(e)\pi(f)\pi(e)$, so
    \begin{align*}
    \cS_{v,w,f}(\lambda)&=\langle\pi(f)\,x_\lambda,\,y_\lambda\rangle, \text{ where}\\
    x_\lambda&:=\pi(e)\pi(\varpi^\lambda)v\in\pi^{G_{-n}},\qquad
    y_\lambda:=\tilde\pi(e)\tilde\pi(\varpi^\lambda)w\in\tilde\pi^{G_{-n}}.
    \end{align*}
    Since $e$ fixes $e_p^\vee$ and $e_q$, we have $\langle x_\lambda,e_p^\vee\rangle=\langle\pi(\varpi^\lambda)v,e_p^\vee\rangle=m^\pi_{v,e_p^\vee}(\varpi^\lambda)$ and $\langle e_q,y_\lambda\rangle=\langle e_q,\tilde\pi(\varpi^\lambda)w\rangle=m^{\tilde\pi}_{w,e_q}(\varpi^\lambda)$. Expanding $x_\lambda=\sum_p\langle x_\lambda,e_p^\vee\rangle e_p$ and $\pi(f)e_p=\sum_q a_{pq}e_q$, \DimaA{we have}
    $$\cS_{v,w,f}(\lambda)=\sum_{p,q}a_{pq}\,m^\pi_{v,e_p^\vee}(\varpi^\lambda)\,m^{\tilde\pi}_{w,e_q}(\varpi^\lambda).$$
    Finally $m^\pi_{v,e_p^\vee}(\varpi^\lambda)\,m^{\tilde\pi}_{w,e_q}(\varpi^\lambda)=\bar m_{v,e_p^\vee}(\lambda)\,\check m_{w,e_q}(\lambda)$, as the twists $\tilde\chi^{-1}(\lambda)$ and $\tilde\chi(\lambda)$ cancel, giving the assertion with $w_p=e_p^\vee$, $v_q=e_q$.
\end{proof}
\begin{cor}\label{cor:S.adm}
    There is a recursion datum $\mdef{\cR_\cS}$   on $\Lambda^+_{ad}({A_{\bfT}})$  and $\mdef{c_{\cS}}\in \Z_{>0}$ s.t. for any $n\in \Z_{>0}$ and any $v\in \pi^{G_{-n}}$,    $w\in \tilde\pi^{G_{-n}}$, and $f\in \cH_{G_{-n}}(G)$, the element  $\cS_{v,w,f}$ is $(\cR_\cS,n c_{\cS})$-admissible.
\end{cor}
\begin{proof}
    By the previous lemma $\cS_{v,w,f}=\sum_{p,q}a_{pq}\,\bar m_{v,w_p}\cdot\check m_{w,v_q}$, a finite sum. By \Cref{cor:mbar.adm} there are a recursion datum $\cR$ and $c_{asm}\in\Z_{>0}$ (depending only on $\pi$) with each $\bar m_{v,w_p}$ being $(\cR,c_{asm}n)$-admissible. Applying \Cref{cor:mbar.adm} to the irreducible admissible contragredient $\tilde\pi$ (whose twisted matrix coefficient is exactly $\check m$) gives a datum $\tilde\cR$ and $\tilde c_{asm}$ with each $\check m_{w,v_q}$ being $(\tilde\cR,\tilde c_{asm}n)$-admissible. Set $N:=\max(c_{asm},\tilde c_{asm})$, so both factors are, respectively, $(\cR,Nn)$- and $(\tilde\cR,Nn)$-admissible. By \Cref{cor:cone.prod} (applied to $\cR_1=\cR$, $\cR_2=\tilde\cR$) there are a recursion datum $\cR_\cS$ and $c'\in\Z_{>0}$, depending only on $\cR,\tilde\cR$, such that every product $\bar m_{v,w_p}\cdot\check m_{w,v_q}$ is $(\cR_\cS,c'Nn)$-admissible. As the $(\cR_\cS,c'Nn)$-admissible elements of $\C^{\Lambda^+_{ad}({A_{\bfT}})}$ form a $\C$-vector space (admissibility is a family of linear conditions), the finite combination $\cS_{v,w,f}$ is $(\cR_\cS,c_\cS n)$-admissible with $c_\cS:=c'N$.
\end{proof}
\section{Proof of \Cref{thm:recur}}\label{ssec:KY}
\begin{proof}[Idea of the proof]
By
the Cartan decomposition, $\langle\cA_N(m_{v,w}),f\rangle$ is a sum, over the dominant
cone, of the integrals of $m_{v,w}$ over the $G_0$-double cosets.
For simplicity assume that $\Lambda^+(T)=\Lambda^+(A_\bfT)$.

\RamiA{Denote the summand 
attached to a $G_0\times G_0$-coset  corresponding to a coweight $\lambda\in \Lambda^+(T)$ by $\cS^1(\lambda)$. It} is the product of the volume $V(\lambda)$ of a $G_{-n}\times G_{-n}$ coset  inside $G_0\varpi^{\lambda} G_0$ and of the sum
 of matrix coefficients over representatives of 
 $$G_{-n} \backslash (G_0\varpi^{\lambda} G_0)/G_{-n}.$$

  We show that $V(\lambda)$  is admissible (\Cref{lem:V.adm}). We then show that the set $$G_{-n} \backslash (G_0\varpi^{\lambda} G_0)/G_{-n}$$ 
stabilizes along a ray in $\Lambda^+(T)$ after a number of
steps proportional to $n$ (\Cref{lem:stab}).
We then deduce that 
$\cS^1$ is an admissible function (Propositions \ref{prop:adm.int} and \ref{prop:S1.adm}). This allows us to deduce the first part of \Cref{thm:recur} (\EitanD{before the ``Moreover''}).

In general, $\Lambda^+(T)\neq \Lambda^+(A_\bfT)$. Thus we slightly modify the argument above. Namely, we first deal with the cosets represented by $\Lambda^+(A_\bfT)$ as discussed above (in \S \ref{ssec:avg}), and then we use \Cref{lem:fin.index} in order to extend this to the entire $\Lambda^+(T)$ - in \S \ref{ssec:avg.G}.  

For the moreover part of \Cref{thm:recur} we use the
elementary fact that a limit of a sequence satisfying a recursion relation is an effective \RamiC{(fixed)}
finite combination of the elements of the sequence (\Cref{lem:rec}).    
\end{proof}
\subsection{Recursion relation of volumes}\label{ssec:vol}
In this subsection we prove that the volume of a $G_{-n}\times G_{-n}$-coset as a function on $\Lambda^+(A_\bfT)$ is admissible. See \Cref{lem:V.adm} below.

\begin{notation}
    Let $g\in T$ and $n\in \Z_{>0}$. Denote by $\mdef{V_{g,n}}\in \C^{\Lambda^+_{ad}({A_{\bfT}})}$ the element given by 
$$V_{g,n}(\lambda)=\mu_{G^{ad}}(\overline{G_{-n} \varpi^{\tilde \lambda}gG_{-n}})$$
    where $\tilde \lambda\in \Lambda^+({A_{\bfT}})$ is a representative of $\lambda\RamiE{\in \Lambda^+_{ad}({A_{\bfT}})}$ and the bar denotes the image in $G^{ad}$\RamiC{.} 
\end{notation}
\begin{lemma}\label{lem:V.adm}
Let $g\in T^+$. Then 
    there are a recursion datum $\mdef{\cR_g}$ on $\Lambda^+_{ad}({A_{\bfT}})$  and $\mdef{c_{g}}\in \Z_{>0}$ s.t. for any $n\in \Z_{>0}$ the element $V_{g,n}$  is $(\cR_g,n c_{g})$-admissible.
\end{lemma}
\begin{proof}
$ $
\begin{enumerate}[Step 1.]
    \item \RamiA{Construction of a character $\chi_\delta:\Lambda^+_{ad}({A_{\bfT}})\to\C^\times$.}\\
    For a weight $\gamma$ of $\bfA_\bfT$ in $\Lie(\bfG)$ denote by $m_\gamma$ the dimension over $F$ of the corresponding weight space, and set $m_\gamma=0$ if $\gamma$ is not a weight. Let $U^+$ be the unipotent radical of $B$ and let $\delta:T\to \RamiE{\C^\times}$ be the inverse of the modular character of $B$, i.e.
    \begin{equation}\label{eq:delta}    \delta(t)=\bigl|\det\bigl(Ad(t)|_{\Lie(U^+)}\bigr)\bigr|^{-1}=\prod_{\beta\in\Phi_\bfB}\ell^{(m_\beta+2m_{2\beta})\beta(t)},
    \end{equation}    
    the second equality since $\Lie(U^+)=\bigoplus_{\beta\in\Phi_\bfB}\Lie(U_\beta)$ and $\Lie(U_\beta)$ is the sum of the weight spaces of the weights $\beta$ and $2\beta$ by \Cref{not:root.sub}\eqref{not:root.grp}; here $\beta(t)$ is as in \Cref{lem:root.conj}. It is a homomorphism, and it is trivial on $\cZ(G)$ since $Ad$ is trivial on $\cZ(G)$, so $\lambda\mapsto\delta(\varpi^{\tilde\lambda})$ is a well-defined character $\chi_\delta:\Lambda^+_{ad}({A_{\bfT}})\to\C^\times$.
\item \RamiA{$V_{g,n}(\lambda)=\mu_{G^{ad}}(\bar G_{-n})\,\delta(g)\,\chi_\delta(\lambda)$ for every $n\in \Z_{>0}$ and every $\lambda\RamiE{\in \Lambda^+_{ad}({A_{\bfT}})}$.}\\
    Put $t:=\varpi^{\tilde\lambda}g$\RamiE{. As} $g\in T^+$ and $\tilde\lambda\in\Lambda^+(A_\bfT)$ we have $\mathrm{val}\,\alpha(t)=\mathrm{val}\,\alpha(\varpi^{\tilde\lambda})+\mathrm{val}\,\alpha(g)\geq0$ for every $\alpha>0$, i.e.\ $t\in T^+$. 

    \RamiD{Let $U:=\bfU_B(F)$ and $U^-:=\bfU_{B^-}(F)$.}
    By \Cref{lem:iva.decomp} for $P=B$, 
    and 
    taking images in $G^{ad}$ \RamiD{we obtain:}    \begin{equation}\label{lem:V.adm:1}
    \RamiD{
    \bar G_{-n}=\Bigl(\overline{U^{-} \cap G_{-n} }\Bigr)\,\overline{T\cap G_{-n}}\,\Bigl(\overline{U \cap G_{-n} }\Bigr).}
    \end{equation}
    \RamiD{Moreover, an element of the LHS of \eqref{lem:V.adm:1} has exactly one factorization as in the RHS.}

    Conjugation by $\bar t$ fixes $\overline{T\cap G_{-n}}$ \claude{by \Cref{lem:T.norm}}.
    \RamiD{
    By \Cref{lem:U.filtration}\eqref{lem:U.filtration:3} we have
    $$
t\,(U^-\cap G_{-n})\,t^{-1}>U^-\cap G_{-n}
$$ and
    $$
t\,(U\cap G_{-n})\,t^{-1}<U\cap G_{-n}.
$$ 
Thus
$$\bar G_{-n}\cap \bar t\bar G_{-n}\bar t^{-1}=\Bigl(\overline{U^{-}\cap G_{-n}}\Bigr)\,\overline{T\cap G_{-n}}\,\Bigl(\overline{t(U\cap G_{-n}) t^{-1}} \Bigr).$$
    Hence:
    }    
    \begin{align*}
    [\bar G_{-n}:\bar G_{-n}\cap \bar t\bar G_{-n}\bar t^{-1}]
    &\mathrel{\RamiD{=}}\RamiD{[U\cap G_{-n} :t(U\cap G_{-n}) t^{-1}]}\\
    &\mathrel{\RamiD{=}}\RamiD{\mu_{U}(U\cap G_{-n}) /\mu_{U}(t(U\cap G_{-n}) t^{-1})}\\
    &\mathrel{\RamiD{\overset{\eqref{eq:delta}}{=}}}\delta(t),
    \end{align*}
    \RamiD{where $\mu_U$ is a Haar measure on $U$.}
    Finally $\bar G_{-n}\bar t\bar G_{-n}$ is the disjoint union of the cosets $\bar z\bar t\bar G_{-n}$, where $\bar z$ runs over $\bar G_{-n}/(\bar G_{-n}\cap \bar t\bar G_{-n}\bar t^{-1})$, each of measure $\mu_{G^{ad}}(\bar G_{-n})$, so that
    $$V_{g,n}(\lambda)=\mu_{G^{ad}}(\bar G_{-n}\bar t\bar G_{-n})=\mu_{G^{ad}}(\bar G_{-n})\,\delta(t),$$
    and $\delta(t)=\delta(\varpi^{\tilde\lambda})\delta(g)=\chi_\delta(\lambda)\delta(g)$, proving the claim.
    \item \RamiA{Completion of the proof}.\\
    Fix a generating set $\lambda_1,\dots,\lambda_k$ of $\Lambda^+_{ad}({A_{\bfT}})$ and set $$\cR_g:=\bigl((\lambda_1,\dots,\lambda_k),(x-\chi_\delta(\lambda_1),\dots,x-\chi_\delta(\lambda_k))\bigr)$$ and $c_g:=1$. Fix a base $\mu$ and an index $i$. By the \DimaA{previous step,} 
    $$V_{g,n}(\mu+m\lambda_i)=\mu_{G^{ad}}(\bar G_{-n})\delta(g)\chi_\delta(\mu)\chi_\delta(\lambda_i)^m$$ for every $m\in\Z_{\geq0}$, whence $V_{g,n}(\mu+(m+1)\lambda_i)=\chi_\delta(\lambda_i)\,V_{g,n}(\mu+m\lambda_i)$. Thus $m\mapsto V_{g,n}(\mu+m\lambda_i)$ is $(x-\chi_\delta(\lambda_i),c_gn)$-admissible, and as $\mu,i$ were arbitrary $V_{g,n}$ is $(\cR_g,c_gn)$-admissible.
    \end{enumerate}
\end{proof}
\subsection{Stabilization of stabilizers}\label{ssec:stab}
In this subsection we prove that the stabilizer of $G_{-n} \varpi^\lambda G_{-n}$ w.r.t. the action of $G_0\times G_0$ stabilizes along rays in $\Lambda^+(A_\bfT)$. See \Cref{lem:stab} below.

\begin{lemma}\label{lem:N.stab}
Let $t\in T$ and denote by $\bar t\in\Lambda(T)$ its class. Then, with $\cN$, $W$ and $w_\nu$ as in \Cref{not:settings}\eqref{it:Weyl},
$$\cN\cap G_0\cap tG_0t^{-1}=\{\nu\in \cN\cap G_0\ :\ w_\nu(\bar t)=\bar t\}.$$
\end{lemma}
\begin{proof}
Let $\nu\in \cN\cap G_0$ and put $s:=\nu^{-1}t^{-1}\nu t\in T$. Since $\nu^{-1}t^{-1}\nu=w_\nu^{-1}(t^{-1})$, the class of $s$ in $\Lambda(T)$ is $\bar t-w_\nu^{-1}(\bar t)$. We have $t^{-1}\nu t=\nu s$, and $\nu\in G_0$, hence
$$t^{-1}\nu t\in G_0\iff s\in T\cap G_0\overset{\text{Lem }\ref{lem:parahoric}\eqref{lem:parahoric:T}}{=}{}^0T\iff \bar t=w_\nu^{-1}(\bar t),$$
the last equivalence since $\Lambda(T)=T/{}^0T$. Finally $w_\nu^{-1}(\bar t)=\bar t$ iff $w_\nu(\bar t)=\bar t$.
\end{proof}

\begin{lemma}\label{lem:W.stab}
For $z\in \Lambda(T)$ denote $\mdef{W_z}:=\{w\in W: w(z)=z\}$. Let $u,v\in\Lambda^+(T)$ and $k\in\Z_{>0}$. Then $W_{u+kv}=W_u\cap W_v$.
\end{lemma}
\begin{proof}
The inclusion $\supset$ is clear. Since $u,v\in\Lambda^+(T)$, the element $u+kv$ lies in the closed Weyl chamber corresponding to $\bfB$, so $W_{u+kv}$ is generated by the reflections it contains (\cite[Ch.~V, \S3, no.~3, Proposition~2]{Bou}). It is therefore enough to show that $s_\beta\in W_u\cap W_v$ for every $\beta\in\Phi$ with $s_\beta\in W_{u+kv}$. For such a $\beta$, which we may assume to lie in $\Phi_\bfB$, we have $\beta(u+kv)=0$, while $\beta(u)\geq0$ and $\beta(v)\geq0$ by \Cref{not:settings}\eqref{it:M+},\eqref{it:Lam+}; as $k>0$ this forces $\beta(u)=\beta(v)=0$. Since $s_\beta$ acts trivially on $\ker(\beta)$ we get $s_\beta(u)=u$ and $s_\beta(v)=v$.
\end{proof}

\begin{lemma}\label{lem:stab}
    For any $g\in T^+$ and any $\lambda\in \Lambda^+(A_\bfT)$ there is $\mdef{c_{stab,g,\lambda}}\in \Z_{>0}$ s.t. for any $n\in\Z_{>0}$, any $\mu\in \Lambda^+(A_\bfT)$ and any
    $k_1,k_2> n c_{stab,g,\lambda}$, denoting $t_j:=\varpi^{\mu+ k_j\lambda} g$, we have
    \begin{enumerate}
    \item\label{lem:stab:G0} $ G_{-n}(G_0 \cap  t_1 G_0 t_1^{-1})= G_{-n}(G_0 \cap  t_2 G_0 t_2^{-1})$
    \item\label{lem:stab:Gn} $ G_{-n}(G_0 \cap  t_1^{-1} G_{-n} t_1)= G_{-n}(G_0 \cap  t_2^{-1} G_{-n} t_2)$
    \end{enumerate}
\end{lemma}
\begin{proof}
\RamiE{WLOG \claude{assume} $\Phi\neq \emptyset$.}
Let $U^+$ be the unipotent radical of $B$, so that $\Phi=\Phi_{B^-}\sqcup \Phi_B$ by \Cref{not:settings}\eqref{it:Phi.P}. Denote
$$\Phi_1:=\{\beta\in\Phi:\beta(\lambda)>0\},\qquad \Phi_0:=\{\beta\in\Phi:\beta(\lambda)=0\},\qquad \Phi_{-1}:=-\Phi_1 ,$$
so that $\Phi=\Phi_{-1}\sqcup\Phi_0\sqcup\Phi_1$ and, $\lambda$ being dominant, $\Phi_1\subset\Phi_B$ and $\Phi_{-1}\subset\Phi_{B^-}$. Put
$$b:=\max_{\beta\in\Phi}|\beta(g)|,\qquad c_{stab,g,\lambda}:=\RamiD{2}\lceil\fN(b+1)\rceil\in\Z_{>0}.$$

Fix $n\in\Z_{>0}$, $\mu\in\Lambda^+(A_\bfT)$ and $k>nc_{stab,g,\lambda}$, and put $t:=\varpi^{\mu+k\lambda}g\in T^+$, so that
$\beta(t)=\beta(\mu)+k\beta(\lambda)+\beta(g)$ for every $\beta\in\Phi$. By \Cref{not:settings}\eqref{it:N.Theta} we have $\fN\beta(\lambda)\in\Z$, hence $\beta(\lambda)\geq 1/\fN$ for $\beta\in\Phi_1$ and $\beta(\lambda)\leq -1/\fN$ for $\beta\in\Phi_{-1}$; and $\beta(\mu)\geq0$ for $\beta\in\Phi_B$, $\beta(\mu)\leq 0$ for $\beta\in \Phi_{B^-}$. Using $n\geq1$ we get
\begin{equation}\label{lem:stab:1}
\begin{split}
\beta(t)&\geq \frac{k}{\fN}-b>n(b+1)-b\geq n\qquad(\beta\in\Phi_1),\\
\beta(t)&\leq -\frac{k}{\fN}+b<-n<0\qquad(\beta\in\Phi_{-1}),
\end{split}
\end{equation}
\begin{equation}\label{lem:stab:2}
\beta(t)=\beta(\mu)+\beta(g)\quad(\beta\in\Phi_0).
\end{equation}
In particular the right hand side of \eqref{lem:stab:2} does not depend on $k$. We show that both groups in the formulation do not depend on $k$ either.
\begin{itemize}
    \item[\eqref{lem:stab:G0}]
By \Cref{lem:parahoric}\eqref{lem:parahoric:norm}, for any subset $X\subset G_0$ the set $G_{-n}\langle X\rangle$ is a subgroup of $G_0$, and hence
$$G_{-n}\langle X\rangle=\langle G_{-n}\cup X\rangle .$$
Applying this to \Cref{lem:parahoric}\eqref{lem:parahoric:fix} we get
\begin{equation}\label{lem:stab:3}
G_{-n}(G_0\cap tG_0t^{-1})=\bigl\langle G_{-n},\ \cN\cap G_0\cap tG_0t^{-1},\ U_{\beta,\max(0,\beta(t))}\ (\beta\in\Phi)\bigr\rangle .
\end{equation}
For $\beta\in\Phi_1$ we have $\max(0,\beta(t))=\beta(t)\geq n$ by \eqref{lem:stab:1}, hence
$U_{\beta,\max(0,\beta(t))}\subset U_{\beta,n}\subset G_{-n}$ by \Cref{lem:U.filtration}, and these generators may be omitted. For $\beta\in\Phi_{-1}$ we have $\max(0,\beta(t))=0$ by \eqref{lem:stab:1}, and for $\beta\in\Phi_0$ the number $\max(0,\beta(t))$ does not depend on $k$ by \eqref{lem:stab:2}.

Finally, let $\nu_0\in\Lambda^+(T)$ and $\lambda_0\in\Lambda^+(T)$ be the classes of $\varpi^{\mu}g$ and of $\varpi^{\lambda}$. The class of $t$ in $\Lambda(T)$ is $\nu_0+k\lambda_0$, so \Cref{lem:N.stab} and \Cref{lem:W.stab} give
$$\cN\cap G_0\cap tG_0t^{-1}=\{\nu\in\cN\cap G_0:\ w_\nu\in W_{\nu_0}\cap W_{\lambda_0}\},$$
which does not depend on $k$. Thus all the generators in \eqref{lem:stab:3} can be chosen independently of $k$.

\item[\eqref{lem:stab:Gn}]
\RamiD{Let $U:=\bfU_B(F)$ and $U^-:=\bfU_{B^-}(F)$.}
By \Cref{lem:iva.decomp} for $P=B$ 
\RamiD{we have
\begin{equation}\label{lem:stab:iva}
    \RamiD{G_{-n}=(U^-\cap G_{-n})(T\cap G_{-n})(U\cap G_{-n}).}
    \end{equation}
}

\RamiD{
Let $t\in T$ satisfy $\beta(t)\geq0$ for every $\beta\in\Phi_B$.
    By \Cref{lem:U.filtration}\eqref{lem:U.filtration:3} we have
    $$
t^{-1}\,(U^-\cap G_{-n})\,t<U^-\cap G_{-n}<G_{-n}<G_0
$$ 
Thus
\begin{equation}\label{lem:stab:dec}
\begin{split}
G_{0}\cap  t^{-1} G_{-n} t=&\Bigl(t^{-1}(U^{-}\cap G_{-n})t\Bigr)\left({T\cap G_{-n}}\right)\Bigl({t^{-1}(U\cap G_{-n}) t}\cap G_{0} \Bigr)=\\
&\Bigl(t^{-1}(U^{-}\cap G_{-n})t\Bigr)\left({T\cap G_{-n}}\right)\Bigl({t^{-1}(U\cap G_{-n}) t}\cap(U\cap G_{0}) \Bigr).
\end{split}
\end{equation}

Now, we have 
\begin{align*}
 G_{-n}\left(G_{0}\cap  t^{-1} G_{-n} t \right)
&\overset{\eqref{lem:stab:dec}}{=}G_{-n} \Bigl({t^{-1}(U\cap G_{-n}) t}\cap (U \cap G_{0} )   \Bigr)\\
&\overset{\text{Lem }\ref{lem:U.filtration}\eqref{lem:U.filtration:prod}}{=}G_{-n} \Bigl({t^{-1}(\prod_{\beta\in\Phi_B}U_{\beta,n}U_{2\beta,n}) t}\cap (\prod_{\beta\in\Phi_B}U_{\beta,0} )\Bigr)\\
&\overset{\text{Lem }\ref{lem:root.conj}}{=}G_{-n} \Bigl({(\prod_{\beta\in\Phi_B}U_{\beta,n-\beta(t)}U_{2\beta,n-2\beta(t)})}\cap (\prod_{\beta\in\Phi_B}U_{\beta,0} )\Bigr)\\
&\overset{\text{Lem }\ref{lem:U.filtration} \eqref{lem:U.filtration:4}}{=}G_{-n} \Bigl({\prod_{\beta\in\Phi_B}((U_{\beta,n-\beta(t)}U_{2\beta,n-2\beta(t)})\cap U_{\beta,0}})\Bigr)\\
&=G_{-n} \Bigl({\prod_{\beta\in\Phi_1}((U_{\beta,n-\beta(t)}U_{2\beta,n-2\beta(t)})\cap U_{\beta,0}})\Bigr)\\
&\qquad\Bigl({\prod_{\beta\in\Phi_B\cap\Phi_0}((U_{\beta,n-\beta(t)}U_{2\beta,n-2\beta(t)})\cap U_{\beta,0}})\Bigr).
\end{align*}
Now, let $k>n c_{stab,g,\lambda}$ and let 
$t:=\varpi^{\mu+ k\lambda} g$. We \claude{obtain}:
\begin{align*}
G_{-n}\left(G_{0}\cap  t^{-1} G_{-n} t \right)
&=G_{-n} \Bigl({\prod_{\beta\in\Phi_1}((U_{\beta,n-\beta(t)}U_{2\beta,n-2\beta(t)})\cap U_{\beta,0}})\Bigr)\\
&\qquad\Bigl({\prod_{\beta\in\Phi_B\cap\Phi_0}((U_{\beta,n-\beta(t)}U_{2\beta,n-2\beta(t)})\cap U_{\beta,0}})\Bigr)\\
&\overset{\eqref{lem:stab:1},\,\eqref{lem:stab:2}}{=}G_{-n} \Bigl({\prod_{\beta\in\Phi_1}U_{\beta,0}}\Bigr)\\
&\qquad\Bigl({\prod_{\beta\in\Phi_B\cap\Phi_0}((U_{\beta,n-\beta(\varpi^\mu g)}U_{2\beta,n-2\beta(\varpi^\mu g)})\cap U_{\beta,0}})\Bigr).
\end{align*}
This proves the \claude{assertion}.
}
\end{itemize}
\end{proof}

\subsection{Admissibility of the integral of matrix coefficients over $G_0$-double cosets in $G_0 A_\bfT \RamiE{g} G_0$}\label{ssec:avg}

\begin{notation}\label{not:bar}
For a subset $A\subset G$ denote by $\mdef{\bar A}\subset G^{ad}$ its image, and for $z\in G$ denote by $\mdef{\bar z}\in G^{ad}$ its image. Note that $\bar G_0=G^{ad}_0$ and that, by \Cref{lem:parahoric}\eqref{lem:parahoric:norm}, $\bar G_{-n}$ is a normal subgroup of $\bar G_0$.
\end{notation}
\begin{notation}\label{not:S0}
Let
\begin{itemize}
    \item $g\in T^+$,
    \item $\lambda\in\Lambda^+(A_\bfT)$,
    \item $v\in \pi$, $w\in \tilde\pi$,
    \item $f\in \cH(G)$.
\end{itemize}
Denote
    $$\mdef{\cS^0_{v,w,g,f}}(\lambda):=\int_{\overline{G_0\varpi^{\lambda} g G_0}} \langle Ad(h)(m_{v,w}),f\rangle \,\mu_{G^{ad}} .$$
Note that $\cS^0_{v,w,g,f}$ is invariant w.r.t. $\Lambda^0(A_\bfT)$, since $\varpi^{\Lambda^0(A_\bfT)}\subset \cZ(G)$ by \Cref{lem:Lam0.center} and $Ad$ is trivial on $\cZ(G)$; so we will also consider it as an element of $\C^{\Lambda^+_{ad}(A_\bfT)}$. 
\end{notation}

In this subsection we prove
\begin{prop}\label{prop:adm.int}
For any $g\in T^+$
    there are a recursion datum $\mdef{\cR_{\cS_0}(g)}$   on $\Lambda^+_{ad}({A_{\bfT}})$  and $\mdef{c_{\cS_0}(g)}\in \Z_{>0}$ s.t. for any $n\in \Z_{>0}$ and any $v\in \pi^{G_{-n}}$,    $w\in \tilde\pi^{G_{-n}}$, and $f\in \cH_{G_{-n}}(G)$, the element  $\cS^0_{v,w,g,f}$ is $(\cR_{\cS_0}(g),n c_{\cS_0}(g))$-admissible.
\end{prop}
For the proof we will need some preparations:

\begin{lemma}\label{lem:S0.decomp}
Let
\begin{itemize}
    \item $g\in T^+$ and $\lambda\in\Lambda^+(A_\bfT)$; put $x:=\varpi^\lambda g$
    \item $n\in \Z_{>0}$
    \item $v\in \pi^{G_{-n}}$, $w\in \tilde\pi^{G_{-n}}$
    \item $f\in \cH_{G_{-n}}(G)$
\end{itemize}
Choose a set $R\subset \bar G_0$ of representatives of $$\bar G_0/\bar G_{-n}(\bar G_0 \cap  \bar x \bar G_0 \bar x^{-1})$$ and a set $R'\subset \bar G_0$ of representatives of $$\bar G_{-n}(\bar G_0 \cap  \bar x^{-1} \bar G_{-n} \bar x)\backslash \bar G_0 .$$
For $y\in R'$ let $\tilde y\in G_0$ be a lift of $y$. Then we have
$$\cS^0_{v,w,g,f}(\lambda)=V_{g,n}(\lambda)\sum_{h\in R}\sum_{y\in R'} \cS_{g\tilde yv,\,g\tilde yw,\,Ad(h^{-1})(f)}(\lambda).$$
\end{lemma}
\begin{proof}
The summands do not depend on the choice of the lifts $\tilde y$, since $\cZ(G)$ acts on $\pi$ by $\chi$ and on $\tilde\pi$ by $\chi^{-1}$. Put
\begin{align*}
S:=\bar G_0 \cap  \bar x \bar G_0 \bar x^{-1},\qquad S'&:=\bar G_0 \cap  \bar x^{-1} \bar G_{-n} \bar x,\\
F(z)&:=\langle Ad(z)(m_{v,w}),f\rangle\quad (z\in G^{ad}) .
\end{align*}
By \Cref{not:bar} the sets $\bar G_{-n}S$ and $\bar G_{-n}S'$ are subgroups of $\bar G_0$, so that $R$ and $R'$ are well defined. All the computations below take place in $G^{ad}$.
\begin{enumerate}[Step 1.]
    \item\label{lem:S0.decomp:1} 
For $u_1,u_2\in \bar G_{-n}$ and $z\in G^{ad}$ we have $F(u_1zu_2)=F(z)$.\\
Indeed, lifting $u_1,u_2$ to $G_{-n}$, the measure $f$ is invariant under conjugation by $u_1$, being bi-$G_{-n}$-invariant, and $\pi(u_2)v=v$, $\tilde\pi(u_2)w=w$; hence
\begin{align*}
F(u_1zu_2)&=\int_{\RamiE{s\in G^{ad}}} m_{v,w}(u_2^{-1}z^{-1}u_1^{-1}su_1zu_2)\,\RamiE{f}(s)\\
&=\int_{\RamiE{s\in G^{ad}}} m_{v,w}(u_2^{-1}z^{-1}szu_2)\,\RamiE{f}(s)=F(z).
\end{align*}
\RamiE{
\item\label{lem:S0.decomp:1.1} $F(h\bar x y)=\cS_{g\tilde yv,g\tilde yw,Ad(h^{-1})(f)}(\lambda)$ for $h,y\in\bar G_0$.\\
Recall that $Ad(z)(m_{v,w})=m_{zv,zw}$. Since $Ad$ is a homomorphism we have:
\begin{align*}
F(h\bar xy)&=\langle Ad(h)\bigl(Ad(\overline{\varpi^{\lambda}g\tilde y})(m_{v,w})\bigr),f\rangle\\
&=\langle Ad(\varpi^\lambda)(m_{g\tilde yv,g\tilde yw}),Ad(h^{-1})(f)\rangle=\cS_{g\tilde yv,g\tilde yw,Ad(h^{-1})(f)}(\lambda).
\end{align*}
}
\item\label{lem:S0.decomp:2} $\bar G_0\bar x\bar G_0=\bigsqcup_{h\in R}h\bar G_{-n}\bar x\bar G_0$.\\
Since $\bar G_0=\bigsqcup_{h\in R}h\bar G_{-n}S$ and $S\bar x\subset \bar x\bar G_0$ we get
$$\bar G_0\bar x\bar G_0=\bigcup_{h\in R}h\bar G_{-n}S\bar x\bar G_0=\bigcup_{h\in R}h\bar G_{-n}\bar x\bar G_0 .$$
If $h_1\bar G_{-n}\bar x\bar G_0\cap h_2\bar G_{-n}\bar x\bar G_0\neq\emptyset$ then $h_2^{-1}h_1=u_1\bar xz\bar x^{-1}u_2$ for some $u_1,u_2\in \bar G_{-n}$ and $z\in \bar G_0$; since $h_2^{-1}h_1\in \bar G_0$ this gives $\bar xz\bar x^{-1}\in S$ and hence $h_2^{-1}h_1\in \bar G_{-n}S$, i.e. $h_1=h_2$.

\item\label{lem:S0.decomp:3} $\bar G_{-n}\bar x\bar G_0=\bigsqcup_{y\in R'}\bar G_{-n}\bar xy\bar G_{-n}$.\\
Since $\bar G_0=\bigsqcup_{y\in R'}\bar G_{-n}S'y$, $\bar xS'\subset \bar G_{-n}\bar x$, and $y$ normalizes $\bar G_{-n}$, we get
$$\bar G_{-n}\bar x\bar G_0=\bigcup_{y\in R'}\bar G_{-n}\bar xS'\bar G_{-n}y=\bigcup_{y\in R'}\bar G_{-n}\bar xy\bar G_{-n}.$$
If $\bar G_{-n}\bar xy_1\bar G_{-n}=\bar G_{-n}\bar xy_2\bar G_{-n}$, write $\bar xy_1=u_1\bar xy_2u_2$ with $u_1,u_2\in \bar G_{-n}$. Then $y_1u_2^{-1}y_2^{-1}=\bar x^{-1}u_1\bar x$ lies in $S'$, hence $y_1\in S'y_2u_2=S'(y_2u_2y_2^{-1})y_2\subset \bar G_{-n}S'y_2$, i.e. $y_1=y_2$.

\item\label{lem:S0.decomp:4} $\mu_{G^{ad}}(\bar G_{-n}h\bar xy\bar G_{-n})=V_{g,n}(\lambda)$ for all $h,y\in \bar G_0$.\\
Indeed $\mu_{G^{ad}}(\bar G_{-n}z\bar G_{-n})=\mu_{G^{ad}}(\bar G_{-n})[\bar G_{-n}:\bar G_{-n}\cap z^{-1}\bar G_{-n}z]$ for every $z\in G^{ad}$, and for $z=h\bar xy$ we have, $h$ and $y$ normalizing $\bar G_{-n}$,
$$\bar G_{-n}\cap z^{-1}\bar G_{-n}z=y^{-1}\bigl(\bar G_{-n}\cap \bar x^{-1}\bar G_{-n}\bar x\bigr)y ,$$
so that $\mu_{G^{ad}}(\bar G_{-n}z\bar G_{-n})=\mu_{G^{ad}}(\bar G_{-n}\bar x\bar G_{-n})=V_{g,n}(\lambda)$.

\item Completion of the proof.\\
By Steps \ref{lem:S0.decomp:2} and \ref{lem:S0.decomp:3} and the normality of $\bar G_{-n}$ in $\bar G_0$ we obtain the disjoint decomposition
$$\overline{G_0xG_0}=\bar G_0\bar x\bar G_0=\bigsqcup_{h\in R}\bigsqcup_{y\in R'}\bar G_{-n}\,h\bar xy\,\bar G_{-n},$$
whence, by Steps \ref{lem:S0.decomp:1},\ref{lem:S0.decomp:1.1} and \ref{lem:S0.decomp:4},
\begin{align*}
\cS^0_{v,w,g,f}(\lambda)&=\sum_{h\in R}\sum_{y\in R'}F(h\bar xy)\,\mu_{G^{ad}}(\bar G_{-n}h\bar xy\bar G_{-n})\\
&=V_{g,n}(\lambda)\sum_{h\in R}\sum_{y\in R'}\cS_{g\tilde yv,g\tilde yw,Ad(h^{-1})(f)}(\lambda).
\end{align*}
\end{enumerate}

\end{proof}

\begin{lemma}\label{lem:center}
Let $x\in T$ and $n\in \Z_{\geq 0}$. Then
$$\bar G_0\cap \bar x\bar G_0\bar x^{-1}=\overline{G_0\cap xG_0x^{-1}},\qquad
\bar G_0\cap \bar x^{-1}\bar G_{-n}\bar x=\overline{G_0\cap x^{-1}G_{-n}x}.$$
\end{lemma}
\begin{proof}
The inclusions $\supset$ are clear. Let $z\in \bar G_0\cap \bar x\bar G_0\bar x^{-1}$ and write $z=\overline{g_0}=\overline{xg_0'x^{-1}}$ with $g_0,g_0'\in G_0$; then $d:=(xg_0'x^{-1})^{-1}g_0$ lies in $\cZ(G)$. Both $G_0$ and $xG_0x^{-1}$ are compact subgroups of $G$, hence are contained in ${}^0G$ by \Cref{not:settings}\eqref{it:0H}, and so $d\in \cZ(G)\cap{}^0G$. As $\cZ(\bfG)^\circ\to \bfG^{ab}$ is an isogeny, the group $\cZ(G)\cap{}^0G$ is compact; being a compact subgroup of $T$ it satisfies
$$\cZ(G)\cap{}^0G\subset {}^0T\overset{\text{Lem }\ref{lem:parahoric}\eqref{lem:parahoric:T}}{=}T\cap G_0$$
by \Cref{not:settings}\eqref{it:0H}. Hence $g_0d^{-1}=xg_0'x^{-1}\in G_0\cap xG_0x^{-1}$ and $z=\overline{g_0d^{-1}}$. The second equality is proved in the same way, $G_{-n}$ being compact as well.
\end{proof}

\begin{proof}[Proof of \Cref{prop:adm.int}]
Fix $\lambda_1,\dots,\lambda_r\in\Lambda^+(A_\bfT)$ representing a generating set of the monoid $\Lambda^+_{ad}(A_\bfT)$, and let
\begin{itemize}
\item $c_1:=\max_{i}c_{stab,g,\lambda_i}$, with $c_{stab,g,\lambda_i}$ as in \Cref{lem:stab};
\item $c_2\in\Z_{>0}$ be as in \Cref{lem:T.conj} for the finite set $X=\{g,g^{-1}\}$;
\item $q_i$ (resp. $q_i'$) be the polynomial given by \Cref{lem:prop} for the recursion datum $\cR_{\cS}$ of \Cref{cor:S.adm} (resp. $\cR_g$ of \Cref{lem:V.adm}) and the element $\lambda_i$;
\item $p_i\in\C[x]$ and $c_3\in\Z_{>0}$ be given by \Cref{lem:1d.prod} for the pair $q_i,q_i'$.
\end{itemize}
Put $\cR_{\cS_0}(g):=\bigl((\lambda_1,\dots,\lambda_r),(p_1,\dots,p_r)\bigr)$ and $c_{\cS_0}(g):=\max(c_1,c_3c_2c_{\cS},c_3c_g)$.

Fix $n,v,w,f$ as in the formulation, an index $i$ and $\mu\in\Lambda^+(A_\bfT)$, and put $x_k:=\varpi^{\mu+k\lambda_i}g$. We have to show that the function $k\mapsto \cS^0_{v,w,g,f}(\mu+k\lambda_i)$ is $(p_i,c_{\cS_0}(g)n)$-admissible.

By \Cref{lem:center} \DimaA{we have}
\begin{align*}
\bar G_{-n}(\bar G_0\cap \bar x_k\bar G_0\bar x_k^{-1})&=\overline{G_{-n}(G_0\cap x_kG_0x_k^{-1})},\\
\bar G_{-n}(\bar G_0\cap \bar x_k^{-1}\bar G_{-n}\bar x_k)&=\overline{G_{-n}(G_0\cap x_k^{-1}G_{-n}x_k)}.
\end{align*}
\DimaA{By \Cref{lem:stab} these groups do not depend on $k$ as long as $k>c_1n$}. Fix such a $k$ and let $R,R'$ be as in \Cref{lem:S0.decomp} for $x_{k}$; then $R,R'$ satisfy the requirements of \Cref{lem:S0.decomp} for $x_{k'}$ for every $k'>c_1n$, and hence
\begin{equation}\label{prop:adm.int:1}
\cS^0_{v,w,g,f}(\mu+k'\lambda_i)=\cV(k')\cdot\Sigma(k')\qquad(k'>c_1n),
\end{equation}
where $\cV(k'):=V_{g,n}(\mu+k'\lambda_i)$ and $\Sigma(k'):=\sum_{h\in R}\sum_{y\in R'}\cS_{g\tilde yv,g\tilde yw,Ad(h^{-1})(f)}(\mu+k'\lambda_i)$.

For $y\in R'$ we have $\pi(\tilde y)v\in \pi^{G_{-n}}$ and $\tilde\pi(\tilde y)w\in\tilde\pi^{G_{-n}}$ by \Cref{lem:parahoric}\eqref{lem:parahoric:norm}, and $g^{-1}G_{-c_2n}g<G_{-n}$ by \Cref{lem:T.conj}, so that $g\tilde yv\in\pi^{G_{-c_2n}}$ and $g\tilde yw\in\tilde\pi^{G_{-c_2n}}$. Also $Ad(h^{-1})(f)\in\cH_{G_{-n}}(G)\subset \cH_{G_{-c_2n}}(G)$, since $h$ normalizes $\bar G_{-n}$ by \Cref{lem:parahoric}\eqref{lem:parahoric:norm}. Hence, by \Cref{cor:S.adm}, each $\cS_{g\tilde yv,g\tilde yw,Ad(h^{-1})(f)}$ is $(\cR_\cS,c_\cS c_2n)$-admissible, and so is their (finite) sum, admissibility being a family of linear conditions. By the choice of $q_i$ the function $\Sigma$ is therefore $(q_i,c_\cS c_2n)$-admissible, and by \Cref{lem:V.adm} and the choice of $q_i'$ the function $\cV$ is $(q_i',c_gn)$-admissible. Both are thus $(q_i,N)$- resp. $(q_i',N)$-admissible for $N:=\max(c_\cS c_2,c_g)n$, so by the choice of $p_i$ and $c_3$ the function $\cV\Sigma$ is $(p_i,c_3N)$-admissible.

Let now $k'>c_{\cS_0}(g)n$. Then $k'+j>c_1n$ for every $j\geq0$, so by \eqref{prop:adm.int:1} the values of $k\mapsto \cS^0_{v,w,g,f}(\mu+k\lambda_i)$ at the points $k'+j$ coincide with those of $\cV\Sigma$; since $k'>c_3N$, the recursion defined by $p_i$ holds for $\cV\Sigma$ at $k'$, and hence also for $k\mapsto \cS^0_{v,w,g,f}(\mu+k\lambda_i)$. As $\mu$ and $i$ were arbitrary, $\cS^0_{v,w,g,f}$ is $(\cR_{\cS_0}(g),c_{\cS_0}(g)n)$-admissible.
\end{proof}

\subsection{Admissibility of the integral of matrix coefficients over $G_0$-double cosets in $G$}\label{ssec:avg.G}
\RamiE{In this section we introduce the function $\cS^1_{v,w,f}$ on the collection of $\overline G_0$-\claude{double} cosets \claude{in} $G^{ad}$. This function \claude{assigns} to a coset $D$ the integral 
$$\int_{D} \langle Ad(h)(m_{v,w}),f\rangle \,\mu_{G^{ad}}.$$ See \Cref{not:S1}.

\claude{Following} \Cref{not:settings}\eqref{it:ov1} we \claude{identify} \claude{this} collection of cosets with 
${\Lambda^+_{ad}(T)}:=\Lambda^+(T)/        {\Lambda_{\cZ(G)}(T)}.$

We aim to show that $\cS^1_{v,w,f}$ is \claude{admissible}. However ${\Lambda^+_{ad}(T)}$ is \claude{not} a cone since it might have \claude{torsion}. So, we decompose it to \claude{finitely} many cones, and prove that $\cS^1_{v,w,f}$ is \claude{admissible} on each of them. See \Cref{prop:S1.adm}.
} 

\begin{notation}\label{not:T.ad}$ $
    \begin{enumerate}        \item\label{it:T.ad} \RamiE{Denote $\mdef{\Lambda_{ad}(T)}:=\Lambda(T)/        {\Lambda_{\cZ(G)}(T)}
        $ and} 
        $\mdef{\Lambda^+_{ad}(T)}:=\Lambda^+(T)/        \RamiE{\Lambda_{\cZ(G)}(T)}
        $, \RamiE{following}  \Cref{not:settings}\eqref{it:ZG}. For $x\in \Lambda^+_{ad}(T)$ denote by $\mdef{D_x}\subset G^{ad}$ the $G^{ad}_0\times G^{ad}_0$-double coset that corresponds to $x$ under the bijection of \Cref{not:settings}\eqref{it:ov1}.
        \item\label{it:T.ad0} Denote by $\mdef{\Lambda^0_{ad}(T)}$ the group of invertible elements of the monoid $\Lambda^+_{ad}(T)$. 
        \item\label{it:T.adp} Denote $\mdef{\Lambda'_{ad}(T)}:=\Lambda^+_{ad}(T)/\Lambda^0_{ad}(T)$ and denote by $\mdef{\varrho}:\Lambda^+_{ad}(T)\to \Lambda'_{ad}(T)$ the quotient map.
        \item\claude{\label{it:T.ad.iota}} \RamiE{Denote by $\mdef{\iota}: \Lambda_{ad}(A_\bfT)\to \Lambda_{ad}(T)$ the natural map. It is well defined by \Cref{lem:Lam0.center}}
    \end{enumerate}
\end{notation}

\begin{lemma}\label{lem:T.split}$ $
    \begin{enumerate}
        \item\label{lem:T.split:fin} $\Lambda^0_{ad}(T)$ is the image in $\Lambda^+_{ad}(T)$ of $$\mdef{\Lambda^0(T)}:=\{\lambda\in \Lambda(T)\,|\,\forall \beta\in \Phi \text{ we have }\beta(\lambda)=0\},$$ and $\Lambda^0_{ad}(T)$ is finite. Denote by $\mdef{e_{ad}}\in \Z_{>0}$ the exponent of $\Lambda^0_{ad}(T)$.
        \item\label{lem:T.split:cone} $\Lambda(T)/\Lambda^0(T)\RamiE{\cong \Lambda_{ad}(T)/\Lambda_{ad}^0(T)}$ is a lattice 
        \RamiE{and $\Lambda'_{ad}(T)\subset \Lambda_{ad}(T)/\Lambda_{ad}^0(T)$}
        is a cone whose only invertible element is $0$ and which generates \RamiE{$\Lambda_{ad}(T)/\Lambda_{ad}^0(T)$}. 
        \item\label{lem:T.split:sec} There is a homomorphism of monoids $\mdef{\sigma}:\Lambda'_{ad}(T)\to \Lambda^+_{ad}(T)$ with $\varrho\circ \sigma=\id$, and 
        $$\Lambda^+_{ad}(T)=\bigsqcup_{q\in \Lambda^0_{ad}(T)}\bigl(q+\sigma(\Lambda'_{ad}(T))\bigr).$$
        \item\label{lem:T.split:Theta} The function $\Theta$ of \Cref{not:settings}\eqref{it:N.Theta} vanishes on $\Lambda^0_{ad}(T)$. The induced function $\mdef{\Theta'}:\Lambda'_{ad}(T)\to \Z_{\geq 0}$ is a homomorphism of monoids with trivial kernel, and $\Theta(x)=\Theta'(\varrho(x))$ for every $x\in \Lambda^+_{ad}(T)$.
        \item\label{lem:T.split:A} 
        $\iota$ is an \claude{embedding} and its image is of \claude{finite} index in 
        $\Lambda_{ad}(T)$.
        
    \end{enumerate}
\end{lemma}
\RamiE{
\begin{proof}
    $ $
    \begin{enumerate}[(1)]
        \item The only non trivial part is the \claude{finiteness} of $\Lambda_{ad}^0(T)$. So we need to show that $        {\Lambda_{\cZ(G)}(T)}$ has \claude{finite} index in $\Lambda^0(T)$. By \Cref{lem:Lam0.center} we have ${\Lambda_{\cZ(G)}(T)}\supset \Lambda_{0}(A_\bfT)$. The \claude{statement} now follows from  \cite[V.3.17]{Ren}.
        \item Clear.
        \item We can decompose $\Lambda_{ad}(T)= \Lambda^0_{ad}(T)\oplus \Gamma$ \claude{where} $\Gamma$ is free. \claude{Extend} $\varrho$ to $\tilde\varrho: \Lambda_{ad}(T)\to \Lambda_{ad}(T)/\Lambda^0_{ad}(T)$. \claude{The} map  $\tilde\varrho$ \claude{induces} an \claude{isomorphism} $\Gamma\to \Lambda_{ad}(T)/\Lambda^0_{ad}(T)$. Let $\tilde\sigma$ be its \claude{inverse}. \claude{Taking} $\sigma:=\tilde\sigma|_{\Lambda'_{ad}(T)}$ \claude{completes} the proof.
        \item Clear.
        \item Follows from  \cite[V.3.17]{Ren}.
    \end{enumerate}
\end{proof}
}
\begin{notation}\label{not:S1}
    Let $v\in \pi$, $w\in \tilde\pi$ and $f\in \cH(G)$. Denote by $\mdef{\cS^1_{v,w,f}}\in \C^{\Lambda^+_{ad}(T)}$ the function given by
    $$\cS^1_{v,w,f}(x):=\int_{D_x} \langle Ad(h)(m_{v,w}),f\rangle \,\mu_{G^{ad}} ,$$
    with $D_x$ as in \Cref{not:T.ad}\eqref{it:T.ad}.
\end{notation}

\begin{lemma}\label{lem:S1.S0}
    For any $g\in T^+$, $\lambda\in \Lambda^+(A_{\bfT})$, $v\in \pi$, $w\in \tilde\pi$ and $f\in \cH(G)$ we have
    $$\cS^1_{v,w,f}\bigl([g]+\iota([\lambda])\bigr)=\cS^0_{v,w,g,f}(\lambda),$$
    where $[g]\in \Lambda^+_{ad}(T)$ and $[\lambda]\in \Lambda^+_{ad}(A_{\bfT})$ denote the classes of $g$ and of $\lambda$, and $\iota$ is as in \claude{\Cref{not:T.ad}\eqref{it:T.ad.iota}}.
\end{lemma}
\begin{proof}
Put $t:=\varpi^\lambda g$. For every $\beta\in \Phi_\bfB$ we have $\beta(t)=\beta(\lambda)+\beta(g)\geq 0$, so $t\in T^+$. The class of $t$ in $\Lambda(T)$ is the sum of the classes of $\varpi^\lambda$ and of $g$, so its class in $\Lambda^+_{ad}(T)$ is $[g]+\iota([\lambda])$. Hence, by \Cref{not:settings}\eqref{it:cartan},\eqref{it:ov1} and \Cref{not:T.ad}\eqref{it:T.ad},
$$D_{[g]+\iota([\lambda])}=\overline{G_0\varpi^{\lambda}gG_0} ,$$
and the two sides of the asserted equality are integrals of the same function over the same set.
\end{proof}

\begin{cor}\label{cor:S1.coset}
    Let $g\in T^+$ and let $\cR_{\cS_0}(g)$ and $c_{\cS_0}(g)$ be as in \Cref{prop:adm.int}. Then for any $n\in \Z_{>0}$, $v\in \pi^{G_{-n}}$, $w\in \tilde\pi^{G_{-n}}$ and $f\in \cH_{G_{-n}}(G)$, the function $\Lambda^+_{ad}(A_{\bfT})\to \C$ given by $\lambda\mapsto \cS^1_{v,w,f}([g]+\iota(\lambda))$ is $(\cR_{\cS_0}(g),c_{\cS_0}(g)n)$-admissible.
\end{cor}
\begin{proof}
By \Cref{lem:S1.S0} the function in question sends $\lambda$ to $\cS^0_{v,w,g,f}(\tilde \lambda)$, where $\tilde\lambda\in \Lambda^+(A_\bfT)$ is any representative of $\lambda$; that is, it is the element $\cS^0_{v,w,g,f}$ of $\C^{\Lambda^+_{ad}(A_\bfT)}$ of \Cref{not:S0}. The assertion is thus \Cref{prop:adm.int}.
\end{proof}

\begin{prop}\label{prop:S1.adm}
    There are a recursion datum $\mdef{\cR_{\cS_1}}$ on $\Lambda'_{ad}(T)$ and $\mdef{c_{\cS_1}}\in \Z_{>0}$ s.t. for any $n\in \Z_{>0}$, any $v\in \pi^{G_{-n}}$, $w\in \tilde\pi^{G_{-n}}$, $f\in \cH_{G_{-n}}(G)$ and any $q\in \Lambda^0_{ad}(T)$, the function $\Lambda'_{ad}(T)\to \C$ given by $\lambda\mapsto \cS^1_{v,w,f}(q+\sigma(\lambda))$ is $(\cR_{\cS_1},c_{\cS_1}n)$-admissible.
\end{prop}
\begin{proof}[Idea of the proof]
\RamiF{
Consider the diagram}
$$
\begin{tikzcd}[row sep=large, column sep=large]
e_{ad}\Lambda^+_{ad}(A_\bfT)\arrow[d,hook,"\varrho \circ \iota"]\arrow[dr,hook,"\iota"] & \\
\Lambda'_{ad}(T)\arrow[r,hook,"\sigma"] & \Lambda^+_{ad}(T)
\end{tikzcd}
$$
\RamiF{
It is easy to see that this diagram is commutative. Note that the multiplication by $e_{ad}$ is required to guarantee the commutativity of the diagram.

We need to show that a certain function on $\Lambda'_{ad}(T)$ is admissible. So, by \Cref{lem:fin.index} it is enough to show that its restriction to each coset of $\varrho\circ\iota(e_{ad}\Lambda^+_{ad}(A_\bfT))$ is admissible. This follows from \Cref{cor:S1.coset}.}     
\end{proof}

\begin{proof}
\RamiE{
\claude{\Cref{lem:T.split}\eqref{lem:T.split:A}} implies that there is a \claude{finite} set $Z\subset \Lambda_{ad}'(T)$ s.t. $$Z+\varrho \circ \iota (e_{ad}\Lambda_{ad}^+(A_\bfT))=\Lambda_{ad}'(T).$$

For any $q\in \claude{\Lambda^0_{ad}(T)}$ and $z\in Z$ \claude{choose} $g_{z,q}\in T^+$ s.t. $[g_{z,q}] =q + \sigma(z)$.

Consider the finite collection of recursion data $\{\cR_{\cS_0}(g_{z,q})\}_{z\in Z,q\in \Lambda^0_{ad}(T)}$ on $\Lambda^+_{ad}({A_{\bfT}})$. 
For each of them consider 
the recursion datum $\cR_{\cS_0}(g_{z,q})^{e_{ad}}$ on $e_{ad} \Lambda_{ad}^+(A_\bfT)$ guaranteed by \claude{\Cref{lem:scale.data}}.

Let $\cR_0$, $c_0$ be the recursion datum and integer \claude{guaranteed by \Cref{lem:many.data}} for the collection 
$\{\cR_{\cS_0}(g_{z,q})^{e_{ad}}\}$.

\claude{Applying \Cref{lem:fin.index}} for the pair $\Span_\Z (\varrho \circ \iota (e_{ad}\Lambda_{ad}^+(A_\bfT))) \subset \Span_\Z (\Lambda'_{ad}({{T}}))$ and the recursion datum $\cR_0$ we obtain a recursion datum $\cR$ on $\Lambda'_{ad}({{T}})$ and an integer $c\in \Z_{>0}$.

Now we take $\cR_{\cS\claude{_1}}:=\cR$ and $c_{\cS\claude{_1}}:=c c_0 \max_{q,z} c_{\cS_0}(g_{z,q})$.

Fix $n,v,w,f,q$ as in the proposition. Let
 $\alpha:\Lambda'_{ad}(T)\to \C$ be given by $$\alpha(\lambda):= \cS^1_{v,w,f}(q+\sigma(\lambda)).$$
 For any $z\in Z$ let
 $\beta_z:e_{ad}\Lambda^+_{ad}(A_\bfT)\to \C$ be given by $$\beta_z(\lambda):= \alpha(z+\varrho\circ \iota(\lambda))= \cS^1_{v,w,f}([g_{z,q}]+\iota(\lambda)).$$
 Thus, by \claude{\Cref{cor:S1.coset}} and \claude{\Cref{lem:scale.data}} $\beta_z$ is $(\cR_{\cS_0}(g_{z,q})^{e_{ad}},c_{\cS_0}(g_{z,q})n)$-admissible.
Therefore, by \claude{\Cref{lem:many.data}} it is
$(\cR_{0},c_0\max_{z,q} c_{\cS_0}(g_{z,q})n)$-admissible. Thus, by \claude{\Cref{lem:fin.index}}
$\alpha$ is $(\cR,c c_0\max_{z,q} c_{\cS_0}(g_{z,q})n)$-admissible as \claude{required}.}
\end{proof}

\subsection{Limits of functions with recursion relation}\label{ssec:lim}
In this subsection we record the
elementary fact that a limit of a sequence satisfying a recursion relation is an effective
finite combination of the elements of the sequence (\Cref{lem:rec}). This will be used to deduce the moreover part of \Cref{thm:recur} from its first part together with  the Kazhdan--Yom-Din theorem (\Cref{thm:KY}).

\begin{lemma}\label{lem:rec}
    For any tuple $\alpha_1,\dots,\alpha_r$, not all $0$,  and any $d\in \Z$ there exists a tuple $\beta_1,\dots,\beta_s$ s.t. for any $n\in \Z_{>0}$ and any function $f:\Z_{>0}\to \C$ satisfying:
    \begin{itemize}
        \item $\forall k>n$ we have $\sum_{i=1}^r  \alpha_i f(k+i)=0$
        \item the limit $\lim_{k\to \infty} \frac{f(k)}{k^d}$  exists,
    \end{itemize}
    we have
    $$\sum_{i=1}^s  \beta_i f(k+i)=\lim_{k\to \infty} \frac{f(k)}{k^d}\qquad\text{whenever }k>n+1.$$
\end{lemma}

\begin{proof}
\RamiA{Let $S$
be the shift \RamiA{operator on $\C^{Z_{\geq 0}}$ given by} $(Sh)(k)=h(k+1).$}
\begin{enumerate}[Step 1.]
    \item \RamiA{Construction of a polynomial $\gamma$ s.t. $g:=\gamma(S)f$ is  eventually polynomial.}\\
    Discarding trailing zeros we may assume $\alpha_r\neq0$\DimaA{.}
    Let    $$A(x)=\sum_{i=1}^r\alpha_ix^{i-1}=\alpha_r\prod_j(x-\rho_j)^{m_j}$$
    be the characteristic polynomial of the recursion and let $e\geq0$ be the multiplicity of $0$ as a root of $A$, i.e.\ $\alpha_1=\dots=\alpha_{e}=0$ and $\alpha_{e+1}\neq0$. 
    \RamiA{WLOG assume that $\rho_1=1$, if $A(1)\neq 0$ we still assume $\rho_1=1$ but $m_1=0$.}
    \RamiA{Put
    $$\gamma(x):=\prod_{\rho_j\neq1}(x-\rho_j)^{m_j},\qquad D:=\gamma(S),\qquad g:=Df .$$}    
    \item \RamiA{Construction of a polynomial $Q$ s.t. $g(k)=Q(k)$ for $k\geq n+2$.}\\    
The assumption on $f$ reads $(A(S)f)(m)=0$ for $m\geq n+2$. Writing $A(x)=x^{e}B(x)$ with $B(0)\neq0$, it is equivalent to $(B(S)f)(m)=0$ for $m\geq n+2+e$. The roots of $B$ are the non-zero $\rho_j$ with the same multiplicities, and both its constant term and its leading coefficient are non-zero, so \RamiA{one can write}
    \begin{equation}\label{lem:rec:1} f(k)=\sum_{\rho_j\neq0}P_j(k)\rho_j^k\quad (k\geq n+2+e)
    \end{equation}
    \RamiA{where $P_j$ are polynomials satisfying $\deg P_j<m_j$.}

We have $D=S^{e}\prod_{\rho_j\neq0,1}(S-\rho_j)^{m_j}$, so that $(Df)(k)$ involves the values of $f$ only at the indices $\geq k+e$. For every $j$ with $\rho_j\neq1$ the factor $(S-\rho_j)^{m_j}$ annihilates $k\mapsto P_j(k)\rho_j^k$, since $\deg P_j<m_j$. Hence, by \eqref{lem:rec:1},
    \begin{equation}\label{lem:rec:2}
    g(k)=Q(k)\quad(k\geq n+2),\qquad Q:=D(P_1) .
    \end{equation}
    \item \RamiA{$\deg Q\leq d$ and the coefficient of $k^d$ is $\gamma(1)\lim_{k\to \infty} \frac{f(k)}{k^d}$.}\\ 
    Write $L:=\lim_{k\to\infty}f(k)/k^d$. Expanding $\prod_{\rho_j\neq1}(x-\rho_j)^{m_j}=\sum_tc_tx^t$ we get $g=\sum_tc_tf(\cdot+t)$, and $\sum_tc_t=\gamma\RamiD{(1)}$, this being the value of that polynomial at $x=1$. For every fixed $t$ we have $(k+t)^d/k^d\to1$ and hence $f(k+t)/k^d\to L$, so
    $$\lim_{k\to\infty}\frac{g(k)}{k^d}=\Bigl(\sum_tc_t\Bigr)L=\gamma\RamiD{(1)} L .$$
    By \eqref{lem:rec:2} the left hand side equals $\lim_{k\to\infty}Q(k)/k^d$ for the polynomial $Q$. If $d<0$ this forces $Q=0$ and $\gamma\RamiD{(1)} L=0$, hence $L=0$. If $d\geq0$ it forces $\deg Q\leq d$, and then the coefficient of $k^d$ in $Q$ equals $\gamma\RamiD{(1)} L$.
    \item \RamiA{Completion of the proof.}\\
    Assume first that $d\geq0$ and put
    $$\beta(x):=\frac{x\,(x-1)^d\prod_{\rho_j\neq1}(x-\rho_j)^{m_j}}{d!\,\gamma\RamiD{(1)}}=:\sum_{i=1}^s\beta_ix^i ,$$
    which has no constant term ($\beta(0)=0$) and whose coefficients depend only on $(\alpha_1,\dots,\alpha_r)$ and $d$. Let $k>n+1$. Then $(S-1)^d$, applied at the point $k+1$, involves $g$ only at the indices $\geq k+1\geq n+2$, where $g=Q$ by \eqref{lem:rec:2}. As $(S-1)^d$ is the $d$-th finite difference and $\deg Q\leq d$ by Step 2, $(S-1)^dQ$ is the constant $d!$ times the coefficient of $k^d$ in $Q$, that is $d!\,\gamma\RamiD{(1)} L$. Hence
    \begin{align*}
    \sum_{i=1}^s\beta_i f(k+i)&=\bigl(\beta(S)f\bigr)(k)=\frac{\bigl((S-1)^dg\bigr)(k+1)}{d!\,\gamma\RamiD{(1)}}=\frac{d!\,\gamma\RamiD{(1)} L}{d!\,\gamma\RamiD{(1)}}=L .
    \end{align*}
    If $d<0$, take for $\beta$ the polynomial just defined with $0$ in place of $d$. By Step 2 we have $Q=0$ and $L=0$, so the same computation gives, for $k>n+1$,
    $$\sum_{i=1}^s\beta_i f(k+i)=\frac{g(k+1)}{\gamma\RamiD{(1)}}=\frac{Q(k+1)}{\gamma\RamiD{(1)}}=0=L .$$
 \end{enumerate}

\end{proof}

\subsection{Proof of \Cref{thm:recur}}\label{ssec:pf.recur}
\begin{proof}[Proof of \Cref{thm:recur}]
Let $\cR_{\cS_1}$ and $c_{\cS_1}$ be as in \Cref{prop:S1.adm}, and let $\sigma$ and $\Theta'$ be as in \Cref{lem:T.split}. By \Cref{lem:T.split}\eqref{lem:T.split:Theta} the map $\Theta':\Lambda'_{ad}(T)\to \Z_{\geq0}$ is a homomorphism of monoids with trivial kernel. So \Cref{cor:acc_pushforward}, applied to $\phi=\Theta'$ and $\cR=\cR_{\cS_1}$, provides a nonzero polynomial $p$ and a constant $c'\in\Z_{>0}$ such that, writing
$$p=:\sum_{l=1}^r\alpha_lx^{l-1},\qquad c_{\cF}:=c'c_{\cS_1},\qquad \RamiD{c}:=c_{\cF}+1 ,$$
for any $n\in\Z_{>0}$, $v\in\pi^{G_{-n}}$, $w\in\tilde\pi^{G_{-n}}$, $f\in\cH_{G_{-n}}(G)$ and $q\in \Lambda^0_{ad}(T)$ (for which $\lambda\mapsto \cS^1_{v,w,f}(q+\sigma(\lambda))$ is $(\cR_{\cS_1},c_{\cS_1}n)$-admissible by \Cref{prop:S1.adm}) the function
\begin{equation}\label{thm:recur:2}
\cF_{q}(k):=\sum_{\substack{\lambda\in \Lambda'_{ad}(T)\\ \Theta'(\lambda)\leq k}}\cS^1_{v,w,f}\bigl(q+\sigma(\lambda)\bigr)
\end{equation}
is $(p,c_{\cF}n)$-admissible. (The function $\cF_q$ depends also on $n,v,w$ and $f$, which we suppress from the notation.) In particular $\alpha_1,\dots,\alpha_r$ are not all $0$.

Fix now
    \begin{itemize}
        \item $n\in \Z_{>0}$,
    \item $v \in \pi^{G_{-n}}$,  $w\in \tilde\pi^{G_{-n}}$,
    \item $f\in \cH_{G_{-n}}(G)$,
    \end{itemize}
and let $N\in \Z_{>0}$. By \Cref{not:T.ad}\eqref{it:T.ad}, \Cref{not:settings}\eqref{it:Gn} and \Cref{lem:T.split}\eqref{lem:T.split:sec},\eqref{lem:T.split:Theta},
$$G^{ad}_N=\bigsqcup_{q\in \Lambda^0_{ad}(T)}\ \bigsqcup_{\substack{\lambda\in \Lambda'_{ad}(T)\\ \Theta'(\lambda)\leq N}} D_{q+\sigma(\lambda)} .$$
Hence, by \Cref{not:S1} and \eqref{thm:recur:2},
\begin{equation}\label{thm:recur:3}
\langle\cA_N(m_{v,w}),f\rangle=\int_{G^{ad}_N}\langle Ad(z)(m_{v,w}),f\rangle\,\mu_{G^{ad}}=\sum_{q\in \Lambda^0_{ad}(T)}\cF_q(N).
\end{equation}
Therefore, for $k> c_{\cF} n$,
\begin{align*}
\phi(k)&=\sum_{l=1}^r\alpha_l\,\langle\cA_{k+l}(m_{v,w}),f\rangle
\overset{\eqref{thm:recur:3}}{=}\sum_{q\in \Lambda^0_{ad}(T)}\ \sum_{l=1}^r\alpha_l\,\cF_q(k+l)=0 ,
\end{align*}
the last equality being the recursion defined by $p$ for the functions $\cF_q$, applied at the point $k+1>c_{\cF}n$. In particular $\phi$ is constant (equal to $0$) for $k>\RamiD{c} n$. This proves the first assertion, with $\alpha_1,\dots,\alpha_r$ and $\RamiD{c}$ fixed in advance.

For the ``moreover'' part assume $\pi$ is tempered. By \Cref{thm:KY}, there are $d\in\Z_{\geq 0}$ and $e\in\R_{>0}$ such that, in the space of generalized functions,
\begin{align*}
\langle v,w\rangle \chi_\pi&=e\lim_{n\to\infty}\frac{1}{n^{d}}\int_{G^{ad}_n}(Ad(g)m_{v,w})\,\mu_{G^{ad}}(g),\\
\text{ with } (Ad(g)m)(x)&=m(g^{-1}xg).
\end{align*}
Put $C=1/e$. We get
\begin{equation}\label{=Aklim}  
\lim_{k\to\infty}\frac{\langle\cA_k(m_{v,w}),f\rangle}{k^{d}}=C\,\langle\chi_\pi,f\rangle\,\langle v,w\rangle\qquad\text{for every }n,v,w,f.
\end{equation}
Apply \Cref{lem:rec} to the tuple $(\alpha_1,\dots,\alpha_r)$ and to $d$: it yields $\beta_1,\dots,\beta_s\in\C$, depending only on $\pi$; set $\alpha_i':=\beta_i/C$ (depending only on $\pi$). Now fix any $n,v,w,f$ and put $\psi(k):=\langle\cA_k(m_{v,w}),f\rangle$; by the first part $\psi$ satisfies $\sum_{l=1}^r\alpha_l\psi(k+l)=0$ for $k>c_{\cF}n$, and $\psi$ has the limit above, so \Cref{lem:rec} gives \begin{equation}\label{=betlim}
\sum_{i=1}^s\beta_i\,\psi(k+i)=\lim_{k\to\infty}\psi(k)/k^{d}\end{equation} for $k>c_{\cF}n+1$, in particular for $k>\RamiD{c}n$. Hence, for $k>\RamiD{c}n,$
\begin{align*}
\phi(k)=\sum_{i=1}^s\alpha_i'\,\langle\cA_{k+i}(m_{v,w}),f\rangle\overset{\DimaA{\eqref{=betlim}}}{=}\frac1C\sum_{i=1}^s\beta_i\,\psi(k+i)
\overset{\DimaA{\eqref{=Aklim}}}{=}\langle\chi_\pi,f\rangle\,\langle v,w\rangle,
\end{align*}
which is the ``moreover'' assertion\DimaA{.}
\end{proof}
\appendix
\section{Kazhdan -- Yom-Din Theorem with modified balls}\label{sec:KY}
\Cref{thm:KY}, which we attributed to \cite{KYtempered}, is not identical with
\cite[Proposition~1.10]{KYtempered}.
Namely, there are 4 differences:
\begin{enumerate}
    \item\label{it:diff0} In \cite{KYtempered} they assume that $\pi$  is unitary.
    \item\label{it:diff1} In \cite{KYtempered} they assume that $G$  is semi-simple.
    \item\label{it:diff2} In \cite{KYtempered} the integral is taken over a subset of $G$ instead of $G^{ad}$.
    \item\label{it:diff3} In \cite{KYtempered} the subsets of integration are balls that they denote by $G_{<n}$ which is different from the balls $G^{ad}_{n}$ in \Cref{thm:KY} (even when the group is adjoint). 
\end{enumerate}

These 4 differences are inessential, and the proof in \cite[Proposition~1.10]{KYtempered} is valid also for \Cref{thm:KY}. Let us briefly explain why. 

The difference \eqref{it:diff2} is inessential, since in the semi-simple case, the integral over  $G$ is proportional to the integral over $G^{ad}$. So, we can replace the integration in the formulation of \cite[Proposition~1.10]{KYtempered} with integration of corresponding balls on $G^{ad}$. This makes the formulation meaningful also for the non-semi-simple case.

The difference \eqref{it:diff0} is inessential since any irreducible tempered representation of a reductive group can be twisted by a character resulting in a unitary representation.

If we integrate over subsets of $G^{ad}$ then there is no need to assume that $G$ is semisimple, so difference \eqref{it:diff1} is also inessential.

The justification of \eqref{it:diff3} requires a closer examination of the argument of \cite{KYtempered}.
We claim that the proof in \cite[Proposition~1.10]{KYtempered} is valid for a large family of sequences of balls. To describe this family we introduce some notation:
\begin{notation}
    Let  $$\mdef{\fa}:=(\Lambda(\bfT)\otimes_{\Z}\R)/(\Lambda(\cZ(\bfG))\otimes_{\Z}\R),$$ 
and let $$\mdef{\fa^{+}}:=span_{\R_{\geq 0 }}(\Lambda^+_{ad}(A_\bfT)) \subset \fa.$$
    We can identify $\fa^+$ with $\R_{\geq0}^d$ (cf. \cite[\S4.3]{KYtempered}). This gives us a scalar product on $\fa$.    
    
        Let $\cI\subset \fa^+$ be a finite set s.t. 
        \begin{equation}\label{eq:KY.balls}
        \mdef{\fa^+_{<r}}:=\{x\in\fa^+|\forall v\in \cI \text{ we have } \langle x, v\rangle<r \}
                \end{equation}

is bounded for any $r\in \R_{>0}$ (cf. \cite[Appendix A]{KYtempered}). 

    \Cref{not:settings}\eqref{it:ov1} gives us a map $\mdef{\kappa}: G^{ad}\to \fa^+$.
    Denote $$\mdef{G^{ad,\cI}_{n}}:=\{g\in G^{ad}|\kappa(g)\in \fa^+_{<n}\}.$$    
\end{notation}
\begin{notation}
    We say that 2 sequences of open compact subsets $G^{ad,1}_n,G^{ad,2}_n\subset G^{ad}$  are \tdef{equivalent} if there exists $c\in \Z_{>0}$ satisfying 
$$G^{ad,1}_n\subset G^{ad,2}_{n+c} \subset G^{ad,1}_{n+2c}$$ 
\end{notation}

The proof of \cite[Proposition~1.10]{KYtempered} is valid as a proof of the following statement:
\begin{theorem}\label{thm:KY2}
Let $\cI$ be as above. 
Let $G_n^{ad,*}$ be a sequence of open compact subsets that is equivalent to $G_n^{ad,\cI}$
Let $\pi$ be an irreducible tempered representation of $G$.
Then there are $d\in\Z_{\geq 0}$ and $e\in \R_{>0}$ s.t. for any 
matrix coefficient  $m$ of $\pi$, we have
$$m(1)\chi_\pi=e\lim_{n\to\infty}\frac{1}{n^{d}}\int_{G^{ad,*}_n}(Ad(g)m)\,\mu_{G^{ad}}$$
in the space of generalized functions.
\end{theorem}

\Cref{thm:KY2} generalizes \Cref{thm:KY}. In more detail, we have:
\begin{enumerate}

    \item \label{rem:balls1}
    The sequence of balls $G_{<n}$ in \cite{KYtempered} (after projection to $G^{ad}$) is equivalent to the sequence $G^{ad,\cI_1}_{n}$  for a certain $\cI_1$ as above (see \cite[(4.4)]{KYtempered}).
    \item \label{rem:balls2}
    By definition, our sequence  $G^{ad}_{n}$ is equivalent to the sequence $G^{ad,\cI_2}_{n}$  for a certain $\cI_2$ as above (a singleton)\footnote{In fact it is almost equal to this sequence, the difference is only between $<r$ and $\leq r$ in \eqref{eq:KY.balls}.}.
    \item \label{rem:balls3} However,  $G^{ad,\cI_1}_{n}$  and $G^{ad,\cI_2}_{n}$ are not equivalent.   
\end{enumerate}

\begingroup
  \let\clearpage\relax
  \let\cleardoublepage\relax 
  \printindex
\endgroup

\bibliographystyle{alpha}
\bibliography{Ramibib}
\end{document}